\documentclass{amsart}
\usepackage{amsthm,amsmath,amssymb}
\usepackage{mathrsfs}
\usepackage{booktabs}
\usepackage{caption}
\usepackage{array}
\usepackage[all]{xy}
\newcommand{\ex}[1]{*+<5pt>[o][F]{E_{#1}}}
\newcommand{\li}[1]{*+<3pt>[F]{E_{#1}}}

\newtheorem{theorem}{Theorem}[section]
\newtheorem{lemma}[theorem]{Lemma}

\theoremstyle{definition}
\newtheorem{definition}[theorem]{Definition}

\theoremstyle{definition}

\newtheorem{remark}[theorem]{Remark}

\newtheorem*{xnotation}{Notation}
\newtheorem*{xconjecture}{Conjecture}

\numberwithin{equation}{section}

\begin{document}

\title[The Manin--Peyre conjecture for three del Pezzo surfaces]{The Manin--Peyre conjecture for three del Pezzo surfaces}

\author{Xiaodong Zhao}
\address{Institute of Mathematics, Henan Academy of Sciences, Zhengzhou, Henan 450046, P. R. China}

\email{xdzhao2016@163.com}
\thanks{}


\date{}


\begin{abstract}
The Manin--Peyre conjecture is established for a split singular quintic del Pezzo surface with singularity
type $\mathbf{A}_2$ and two split singular quartic del Pezzo surfaces with singularity types $\mathbf{A}_3+\mathbf{A}_1$ and $\mathbf{A}_4$ respectively.
We use a unified and different slightly method from the previous and improve their error terms. Our method is general and can provide independent proofs of the Manin--Peyre conjecture for most of singular del Pezzo surfaces of degree $d\geq 3$.
\end{abstract}

\maketitle

\section{Introduction}
Let $S\subset \mathbb{P}^d$ be a split singular del Pezzo surface defined over $\mathbb{Q}$ with only $\mathbf{ADE}$-singularities and anticanonically embedded.
And, let $U\subset S$ be the open subset formed by deleting the lines from $S$. In view of the Manin conjecture \cite{FR1}, one predicts
the asymptotic behaviour of the number of rational points of bounded height on $U$, namely
\begin{equation}\label{NUH}
N_{U,H}(B):=\#\{\mathbf{x}\in U(\mathbb{Q}):H(\mathbf{x})\leq B\}= cB(\log B)^{\rho-1}(1+o(1)),
\end{equation}
as $B\rightarrow \infty$, where $\rho$ is the rank of the Picard group of a minimal desingularization of $S$ and $c$ is the constant predicted by Peyre \cite{PE1}. Here $H:\mathbb{P}^d(\mathbb{Q})\rightarrow R_{>0}$ is a height function defined for $(x_0,\ldots,x_d)\in \mathbb{Z}^{d+1}$ such that gcd$(x_0,\ldots,x_d)=1$ by
$$H(x_0,\ldots,x_d)=\max\{|x_0|,\ldots,|x_d|\}.$$

We define the anticanonically embedded del Pezzo surfaces $S_i\subset \mathbb{P}^{d_i}$ with $1\leq i\leq 3$ for $d_1=5$ and $d_2=d_3=4$ over $\mathbb{Q}$ by the following equations:
\begin{equation}\label{S1}
\begin{aligned}
      S_1 : \qquad x_0x_2-x_1x_5 = &x_0x_2-x_3x_4 =
  x_0x_3+x_1^2+x_1x_4\\ = &x_0x_5+x_1x_4+x_4^2 = x_3x_5+x_1x_2+x_2x_4 = 0\quad \text{of type} \;\mathbf{A}_2,
\end{aligned}
\end{equation}
\begin{align}
   S_2 &:& x_0^2+x_0x_3+x_2x_4= x_1x_3- x_2^2= 0 &\quad\text{of type}\;\mathbf{A}_3+\mathbf{A}_1,\label{S2}\\
   S_3 &:& x_0x_1-x_2x_3=x_0x_4+x_1x_2+x_3^2=0 &\quad\text{ of type }\;\mathbf{A}_4.\label{S3}
 \end{align}

The variety $S_1$ has only one singularity, namely $(0 : 0 : 1 : 0 : 0 : 0)$ of type $\mathbf{A}_2$. In addition, it contains precisely four lines, which are defined by
\begin{equation*}
x_0=x_1=x_4=x_i=0 \qquad\text{and}\qquad x_j=x_3=x_5=x_1+x_4=0
\end{equation*}
for $i\in\{3,5\}$ and $j\in\{0,2\}$. Then, $N_{S_1,H}(B)$ will be dominated by the rational points of height at most $B$ that lie on such lines. For this reason, let $U_1$ be the open subset formed by removing the four lines from $S_1$ and $N_{U_1,H}(B)$ be defined as in (\ref{NUH}).

As the minimal desingularisation $\widetilde{S}$ of a split singular del Pezzo surface $S\subset \mathbb{P}^d$ of degree $d\in\{3,\ldots,7\}$ is the blow-up of $\mathbb{P}^2$ in $9-d$ points (see Introduction of \cite{COTS1}), it has Picard group Pic$(\widetilde{S})\cong \mathbb{Z}^{10-d}$. The Cox ring of $\widetilde{S}$ has $13-d$ generators $\eta_1,\ldots,\eta_{13-d}$ with corresponding divisors $E_1,\ldots,E_{13-d}$ (see Theorem 2 and Lemma 12 of \cite{DE1}).

Therefore, if $\widetilde{S}_1$ denotes the minimal desingularization of $S_1$, the Picard group of $\widetilde{S}_1$ has rank $\rho_1=5$.
In addition, we will see in Section 3.1 that the universal torsor is an open subset of the hypersurface
\begin{equation}\label{U01}
\eta_2\eta_5\eta_6^2+\eta_3\eta_7+\eta_4\eta_8=0,
\end{equation}
which is embedded in $\mathbb{A}^8\cong \text{Spec}(\mathbb{Q}[\eta_1,\ldots,\eta_8])$. In 2007, Derenthal \cite{DE3} split their counting into three parts and
used the analytic method to prove the Manin--Peyre conjecture for $S_1$ in (\ref{S1}), i.e.
\begin{equation*}
N_{U_1,H}(B)=c_1B(\log B)^4+O\left(B(\log B)^{4-1/5}\right),
\end{equation*}
where $c_1$ is the Peyre constant.

The variety $S_2$ has has exactly two singularities, namely $(0 : 0 : 0 : 0 : 1)$ of type $\mathbf{A}_3$ and $(0 : 1 : 0 : 0 : 0)$ of type $\mathbf{A}_1$. In addition, it contains precisely three lines, which are defined by
\begin{equation*}
x_0=x_1=x_2=0,\qquad x_0+x_3=x_1=x_2=0,\qquad x_0=x_2=x_3=0.
\end{equation*}
The open subset $U_2$, $\widetilde{S}_2$, $N_{U_2,H}(B)$ are defined in a similar way. Thus, we have $\rho_2=6$. In this case, the universal torsor in \cite{DE2} is an open subset of the hypersurface
\begin{equation} \label{U02}
\eta_1\eta_9+\eta_2\eta_8+\eta_4\eta_5^3\eta_6^2\eta_7=0,
\end{equation}
which is embedded in $\mathbb{A}^9\cong \text{Spec}(\mathbb{Q}[\eta_1,\ldots,\eta_9])$. In 2009, Derenthal \cite{DE2} utilized the universal torsors to obtain the Manin--Peyre conjecture for $S_2$ in (\ref{S2}), i.e.
\begin{equation*}
N_{U_2,H}(B)=c_2B(\log B)^5+O\left(B(\log B)^4(\log\log B)^2\right),
\end{equation*}
where $c_2$ is the Peyre constant.

The variety $S_3$ has only one singularity, namely $(0 : 0 : 0 : 0 : 1)$ of type $\mathbf{A}_4$. In addition, it contains precisely three lines, which are defined by
\begin{equation*}
x_0=x_2=x_3=0,\qquad x_0=x_1=x_3=0,\qquad x_1=x_3=x_4=0.
\end{equation*}
The open subset $U_3$, $\widetilde{S}_3$, $N_{U_3,H}(B)$ are defined in a similar way. Thus, we have $\rho_3=6$. In this case, we will see in Section 5.1 that the universal torsor is an open subset of the hypersurface
\begin{equation} \label{U03}
\eta_5\eta_9+\eta_1\eta_8^2+\eta_3\eta_4^2\eta_6^3\eta_7=0,
\end{equation}
which is embedded in $\mathbb{A}^9\cong \text{Spec}(\mathbb{Q}[\eta_1,\ldots,\eta_9])$. In 2009, Browning and Derenthal \cite{BR1} split their counting into two parts, used the similar method to \cite{DE3} and proved the Manin--Peyre conjecture for $S_3$ in (\ref{S3}), i.e.
\begin{equation*}
N_{U_3,H}(B)=c_3B(\log B)^5+O\left(B(\log B)^{5-2/7}\right),
\end{equation*}
where $c_3$ is the Peyre constant.

In order to to refine the Manin conjecture, Swinnerton-Dyer first formulated a more detailed conjecture.
\begin{xconjecture}(\cite{BR0,SW1})
Let $S'\subset \mathbb{P}^3(\mathbb{Q})$ be a cubic surface. Then there exist positive contstants $\theta_1, \theta_2, \theta_3<1$ with $\theta_1<\min\{\theta_2,\theta_3\}$, a polynomial $f\in\mathbb{R}[x]$ of degree $\rho'-1$, a constant $\gamma\in\mathbb{R}$,
and a sequence of $\gamma_n\in\mathbb{C}$, we have
\begin{equation}\label{N0}
N_{U',H}(B)=Bf(\log B)+\gamma B^{\theta_3}+\mathfrak{R}e\sum{\gamma_nB^{\theta_2+it_n}}+O_\varepsilon(B^{\theta_1+\varepsilon}).
\end{equation}
Here $1/2+it_n$ runs through a set of non-trivial zeros of the Riemann zeta function, with the $t_n>0$ positive and monotonic increasing, such that $\sum{|\gamma_n|}$ and $\sum{t_n^{-2}}$ are convergent.
\end{xconjecture}

\textbf{Conjecture} was originally raised by \cite[Conjecture 2]{SW1} in work of \cite[Conjecture C]{BR0}. Since main term and error term of (\ref{N0}) are more precise, (\ref{N0}) improves (\ref{NUH}). Although this conjecture is related to cubic surfaces, we hope that it is true for any surfaces. Next, we plan to use a unified and different method from \cite{BR1,DE2,DE3} to prove the Manin--Peyre conjecture for $S_1,S_2,S_3$ and improve their error terms with logarithmic saving. Since the better error terms, the results that we obtain lend themselves more readily as a bench test for \textbf{Conjecture}. In fact, \cite[Chapter 5]{BR01} and \cite[Theorem 1.1]{BR02} also utilized the different methods to prove the Manin--Peyre conjecture on del Pezzo surfaces and improve their error terms with logarithmic saving. We can not obtain the error terms with power saving, but our method is general and can generalise to other situations. In \cite{DE1}, Derenthal has classified singular del Pezzo surfaces $S$ of degree $d\geq 3$ where Spec(Cox($\tilde{S}$)) is defined by precisely one torsor equation (see Table \ref{table0} for their singularity types). It is possible to provide independent proofs of the Manin--Peyre conjecture for del Pezzo surfaces on the left in Table \ref{table0} which are called ``Good''. The other surfaces (on the right) are called ``Bad''.  From Table \ref{table0}, we learn that most of surfaces on the left have already been done. However, we can still utilize our method to deal with $\mathbf{A}_3$ degree $5$ del Pezzo surface and $\mathbf{A}_4$ degree $5$ del Pezzo surface to get the better error terms than \cite{CHTS1}.

\begin{table}[!ht]
\caption{Scope of application of the method}
\begin{tabular}{|m{1cm}<{\centering}|m{5.3cm}<{\centering}|m{2cm}<{\centering}|} \hline
Degree & Good & Bad \\ \hline
6& $\mathbf{A}_1$\cite{BR01,CHTS1}, $\mathbf{A}_2$\cite{CHTS1,LO1} & -- \\ \hline
5 & $\mathbf{A}_2$\cite{DE3}, $\mathbf{A}_3$\cite{CHTS1}, $\mathbf{A}_4$\cite{CHTS1,MO1} & $\mathbf{A}_1$\cite{BA3} \\ \hline
4 & $\mathbf{A}_3$\cite{DEFR1}, $\mathbf{A}_3$+$\mathbf{A}_1$\cite{DE2,DEFR2}, $\mathbf{A}_4$\cite{BR1,DEFR2}, 
 $\mathbf{D}_4$\cite{DEFR2,DETS1}, $\mathbf{D}_5$\cite{BRE,CHTS1} & $3\mathbf{A}_1$\cite{LE1}, $\mathbf{A}_2$+$\mathbf{A}_1$\cite{LE1} \\ \hline
3 & $\mathbf{A}_4$+$\mathbf{A}_1$, $\mathbf{A}_5$+$\mathbf{A}_1$\cite{BA1}, $\mathbf{D}_4$\cite{BR03,LE3}, $\mathbf{D}_5$\cite{BR2}, $\mathbf{E}_6$\cite{BRE1} & $\mathbf{A}_3$+$2\mathbf{A}_1$, $2\mathbf{A}_2$+$\mathbf{A}_1$\cite{LE2} \\ \hline
\end{tabular}\\
\label{table0}
\end{table}

In 2009, Derenthal \cite{DE2} developed a new technique for proving asymptotic estimates for the number of integral points of bounded height on universal torsors over del Pezzo surfaces whose counting problems have a special form (see \cite[Section 2]{DE2}). Clearly, $S_1$, $S_2$ and $S_3$ belong the special form. Moreover, the passage to the universal torsor translates the problem of counting rational points on the surface to the problem of counting tuples of integers satisfying the torsor equation and certain height and coprimality conditions. When dealing with the torsor equation, we will encounter one congruence condition. For a linear congruence condition, one can estimate the main term by \cite{DE2}. However, one must show that the error term is negligible, which may or may not be hard. For a quadratic congruence condition, it is more harder.

From our observation, for $S_1$ and $S_2$, the congruence conditions are both linear, since $\eta_7$ occurs linearly in (\ref{U01}) and $\eta_8$ occurs linearly in (\ref{U02}). For $S_3$, the congruence condition is quadratic, since $\eta_8$ occurs with a square power in (\ref{U03}).

\begin{theorem}
For $S_1$, we have
\begin{equation*}
N_{U_1,H}(B)=c_1B(\log B)^4+O\left(B(\log B)^3(\log\log B)\right),
\end{equation*}
where
\begin{equation*}
c_1=\alpha(\widetilde{S}_1)\prod\limits_{p}\left(1-\frac{1}{p}\right)^5\left(1+\frac{5}{p}+\frac{1}{p^2}\right)\omega_\infty(\widetilde{S}_1).
\end{equation*}
Here, both a rational number $\alpha(\widetilde{S}_1)$ and an integral $\omega_\infty(\widetilde{S}_1)$ come from \cite{PE1}, and we will give them explicitly in Section 3.5.
\end{theorem}

\begin{theorem}
For $S_2$, we have
\begin{equation*}
N_{U_2,H}(B)=c_2B(\log B)^5+O\left(B(\log B)^4(\log\log B)\right),
\end{equation*}
where
\begin{equation*}
c_2=\frac{1}{8640}\prod\limits_{p}\left(1-\frac{1}{p}\right)^6\left(1+\frac{6}{p}+\frac{1}{p^2}\right)\omega_\infty(\widetilde{S}_2)
\end{equation*}
and\begin{equation*}
\omega_\infty(\widetilde{S}_2)=\int\limits_{|x_0|,|x_2|,|x_2^2/x_1|,|(x_0^2x_1+x_0x_2^2)/(x_1x_2))|\leq 1,0\leq x_1\leq1}\frac{1}{x_1x_2}\mathrm{d}x_0\mathrm{d}x_1\mathrm{d}x_2.
\end{equation*}
\end{theorem}

\begin{theorem}
For $S_3$, we have
\begin{equation*}
N_{U_3,H}(B)=c_3B(\log B)^5+O\left(B(\log B)^4(\log\log B)\right),
\end{equation*}
where
\begin{equation*}
c_3=\alpha(\widetilde{S}_3)\prod\limits_{p}\left(1-\frac{1}{p}\right)^6\left(1+\frac{6}{p}+\frac{1}{p^2}\right)\omega_\infty(\widetilde{S}_3).
\end{equation*}
Similar to $\alpha(\widetilde{S}_1)$ and $\omega_\infty(\widetilde{S}_1)$, we will give $\alpha(\widetilde{S}_3)$ and $\omega_\infty(\widetilde{S}_3)$ explicitly in Section 5.5.
\end{theorem}

Most proofs of the Manin--Peyre conjecture on del Pezzo surfaces fall into two cases:

(I) Del Pezzo surfaces of some types are toric or equivariant compactifications of the additive group $\mathbb{G}_\text{a}^2$, so the Manin--Peyre conjecture follows from \cite{BA2} or \cite{CHTS1} respectively, by techniques of harmonic analysis on adelic groups.

(II) Without using such a structure, the Manin--Peyre conjecture has been proved on Del Pezzo surfaces remained in some cases via universal torsors.

From \cite[Sections 3.3-3.4]{DE1} and \cite[Theorem]{DELO1}, we learn that $S_i$ with $i=1,2,3$ is neither toric nor equivariant compactifications of $\mathbb{G}_\text{a}^2$. Therefore, $S_i$ do not belong to case (I) and our work is not covered by \cite{BA2,CHTS1}.

As it states in \cite{DE2}, proofs of Manin's conjecture for a split del Pezzo surface $S$ which is defined as above, consist of three main steps.

(1) One constructs an explicit bijection between rational points of bounded height on $S$ and integral points in a region on a universal torsor $\mathcal{T}_{S}$.

(2) Using methods of analytic number theory, one estimates the number of integral points in this region on the torsor by its volume.

(3) One shows that the volume of this region grows asymptotically as predicted by Manin and Peyre.

For $S_1$ and $S_3$, we plan to utilize the idea of \cite{DE2} and the following Lemma \ref{L11} to complete steps (1)-(2) and \cite{DEFR1} to achieve step (3),
 which is different from \cite{BR1} and \cite{DE3}. For $S_2$, the difference from \cite{DE2} is that we use the following Lemma \ref{L11} to deal with the second summation over $\eta_7$. Moreover, for $S_1$ and $S_2$, we utilize a trick from \cite[Lemma 4]{DEFR2} to obtain the sufficiently good error terms. Clearly, Theorems 1.1-1.3 improve the error terms of \cite{BR1,DE2,DE3}.

\begin{xnotation}
The $\mu(n)$, $\phi(n)$, $\tau(n)$ and $\vartheta*\mu(n)$ denote M\"{o}bius function, Euler function, divisor function and Dirichlet convolution separately.
Then, $\phi^\ast(n):=\phi(n)/n$ and $\omega(n)$ denotes the counting functions of the total number of prime factors of $n$, taken without multiplicity.
Moreover, $p^\nu\|n$ means that $\nu\in \mathbb{Z}_{>0}$ such that $p^\nu \mid n$ but $p^{\nu+1} \nmid n$.
\end{xnotation}

\section{Two general preliminary results}
In Section 1, we learn that we can prove the Manin--Peyre conjecture for $S$, which is defined as the above, by steps (1)-(3).
We usually apply the strategy in \cite{DETS1} to complete step (1). For steps (2)-(3), there are many preliminary results in \cite{DE2,DEFR1}.
Here, we will give two general preliminary results for steps (2)-(3) separately, which are useful to deal with the split singular del Pezzo surface over $\mathbb{Q}$.

\subsection{A general preliminary result for step (2)}
$\newline$For $S$, the universal torsor is defined in the variables $\eta_1,\ldots,\eta_{13-d}$. When we deal with the height conditions, by the torsor equation, we may eliminate one variable $\eta_{13-d}$ that occurs linearly in it. Fixing $\eta_1,\ldots,\eta_{12-d}$, in general, we will divide it into three cases to achieve the above step (2):
\begin{itemize}
\item the first summation over $\eta_{12-d}$,
\item the second summation over $\eta_{11-d}$,
\item the remaining summations.
\end{itemize}
We generally utilize \cite[Proposition 2.4]{DE2} and \cite[Proposition 4.3]{DE2} to deal with the first summation and the remaining summations respectively.
Then, the following Lemma \ref{L11} is used to handle the second summation. The proof is similar to \cite[Lemma 3.1]{DE2}, but it may be viewed as the version on $\mathbb{Q}$ of \cite[Proposition 6.1]{DEFR1}. Before proving, we have to quote two definitions from \cite{DE2}.

\begin{definition}(\cite[Definition 6.4]{DE2})
\label{D1}
Let $\Theta_1$ be the set of all non-negative functions $\vartheta: \mathbb{Z}_{>0} \to
  \mathbb{R}$ such that there is a $c\in \mathbb{R}$ and a system of non-negative functions
  $A_p: \mathbb{Z}_{\geq0} \to \mathbb{R}$ for $p \in \mathcal{P}$ satisfying
  \begin{equation*}
    \vartheta(n) = c \prod\limits_{p^\nu \|n}{A_p(\nu)}\prod\limits_{p\nmid n}{A_p(0)}
  \end{equation*}
  for all $n \in \mathbb{Z}$. In this
  situation, we say that $\vartheta \in \Theta_1$ \emph{corresponds to} $c,A_p$.
\end{definition}

\begin{definition}(\cite[Definition 6.6]{DE2})
\label{D2}
 For any $b \in \mathbb{Z}_{>0}$, $C_1,C_2,C_3\in \mathbb{R}_{\ge 1}$, let
  $\Theta_2(b,C_1,C_2,C_3)$ be the set of all functions $\vartheta \in \Theta_1$
  for which there exist corresponding $c, A_p$ satisfying the following
  conditions:
  \begin{enumerate}
  \item For all $p \in \mathcal{P}$ and $\nu \ge 1$,
   $$\left|A_p(\nu)-A_p(\nu-1)\right| \le \begin{cases}
      C_1, &p^\nu \mid b,\\
      C_2p^{-\nu}, &p^\nu \nmid b;
    \end{cases}$$
  \item For all $k \in \mathbb{Z}_{>0}$, we have
    $\left|c\prod\limits_{p\nmid k}A_p(0)\right| \le C_3$.
  \end{enumerate}
  Given $\vartheta\in \Theta_2(b,C_1,C_2,C_3)$, by \cite[Proposition 6.8]{DE2}, for any $q \in \mathbb{Z}_{>0}$, the infinite sum
    and the infinite product
 \begin{equation*}
      \sum_{\substack{n>0\\ \gcd(n,q)}} \frac{(\vartheta*\mu)(n)}{n}, \qquad
      c\prod\limits_{p\nmid q}\left({\left(1-\frac{1}{p}\right) \sum_{\nu=0}^\infty
        \frac{A_p(\nu)}{p^\nu}}\right) \prod\limits_{p\mid q} A_p(0)
    \end{equation*}
  converge to a real number, which we denote as $\mathcal{A}(\vartheta(n),n,q)$.
\end{definition}

If $A_p(\nu)=A_p(\nu+1)$ for all primes $p$ and all $\nu \ge 1$, then the
formula is simplified to
\begin{equation*}
  \mathcal{A}(\vartheta(n),n,q)= c \prod\limits_{p\nmid q}\left((1-\frac{1}{p})A_p(0)+\frac{1}{p}A_p(1)\right) \prod\limits_{p\mid q}{A_p(0)}.
\end{equation*}

\begin{lemma}\label{PR1}
Let $\vartheta \in \Theta_2(b,C_1,C_2,C_3)$ and $q, m\in \mathbb{Z}_{>0}$. There exist $A\in \mathbb{Q}$ and $A_1, A_2 \in \mathbb{Z}$ such that $A=A_1/A_2$, $\gcd(A_2, q)=1$ and $\gcd(\bar{A_1}, q)=1$. Here $\bar{A_1}$ denotes the inverse of $A_1$ modulo $q$. For $t\in\mathbb{R}_{\geq 0}$, we have
\begin{align*}
\sum_{\substack{0<n\leq t}}\vartheta(n)&\sum_{\substack{\varrho^m\equiv An(\bmod q) \\1\leq \varrho\leq q\\ \gcd(\varrho, q)=1}} 1 \\ =&\phi^\ast(q)\mathcal{A}(\vartheta(n),n,q)t
+O\left(\phi(q)\tau(b)(C_1C_2)^{\omega(b)}C_3(\log(t+2))^{C_2}\right).
\end{align*}
\end{lemma}
\begin{proof}
If $\gcd(a, q)=1$, from \cite[Corollary 6.9]{DE2}, we have
\begin{align*}
\sum_{\substack{0<n\leq t\\ n \equiv a(\bmod q)}} \vartheta(n)=\frac{t}{q} \mathcal{A}(\vartheta(n),n,q)+O\left(\tau(b)(C_1C_2)^{\omega(b)}C_3(\log(t+2))^{C_2}\right).
\end{align*}
Changing the order of summation, we obtain
\begin{align*}
\sum_{\substack{0<n\leq t}}\vartheta(n)&\sum_{\substack{\varrho^m\equiv An(\bmod q) \\1\leq \varrho\leq q\\ \gcd(\varrho, q)=1}} 1
=\sum_{\substack{1\leq \varrho\leq q\\ \gcd(\varrho, q)=1}}\sum_{\substack{n\equiv \bar{A}_1A_2\varrho^m(\bmod q) \\0<n\leq t}}\vartheta(n)\\
=&\phi^\ast(q)\mathcal{A}(\vartheta(n),n,q)t+O\left(\phi(q)\tau(b)(C_1C_2)^{\omega(b)}C_3(\log(t+2))^{C_2}\right).
\end{align*}
\end{proof}

\begin{lemma}\label{L11}
Let $q$, $A$, $\varrho$, $m$ and $\vartheta$ be as in the previous lemma, and let $I=[t_1,t_2]$ for $t_1, t_2 \in \mathbb{R}_{\geq 0}$ with $t_1 \le t_2$. Let $g:I \to \mathbb{R}$ be a function that has a continuous derivative whose sign changes only $R(g)$ times on $I$. Then, one has
\begin{align*}
\sum_{\substack{n\in I}} \vartheta(n)\sum_{\substack{\varrho^m\equiv An(\bmod q) \\1\leq \varrho\leq q\\ \gcd(\varrho, q)=1}}g(n)=&\phi^\ast(q)\mathcal{A}(\vartheta(n),n,q)\int\limits_{I}g(t)\mathrm{d}t \\
&+O\left(\phi^\ast(q)\tau(b)(C_1C_2)^{\omega(b)}C_3(\log(|I|+2))^{C_2}M_I(g)\right),
\end{align*}
where
\begin{align*}
M_I(g):=(R(g)+1)\cdot \sup_{t\in I}|g(t)| \qquad\text{and}\qquad
|I|:=\max\left\{|t_1|,|t_2|\right\}.
\end{align*}
\end{lemma}
\begin{proof}
For any $t\in \mathbb{R}_{\geq0}$, let
\begin{align*}
M(t):=\sum_{\substack{0< n\leq t}}\vartheta(n)\sum_{\substack{\varrho^m\equiv An(\bmod q) \\1\leq \varrho\leq q\\ \gcd(\varrho, q)=1}} 1 \quad\text{and}\quad
 S(I):=\sum_{\substack{n\in I}} \vartheta(n)\sum_{\substack{\varrho^m\equiv An(\bmod q) \\1\leq \varrho\leq q\\ \gcd(\varrho, q)=1}}g(n).
\end{align*}
By partial summation, we have
\begin{align*}
S(I)=M(t_2)g(t_2)-M(t_1)g(t_1)-\int\limits_{I}M(t)g'(t)\mathrm{d}t.
\end{align*}
From Lemm \ref{PR1}, we have
\begin{align*}
M(t)=\phi^\ast(q)\mathcal{A}(\vartheta(n),n,q)t+E(t),
\end{align*}
where
\begin{align*}
E(t)=O\left(\phi(q)\tau(b)(C_1C_2)^{\omega(b)}C_3(\log(t+2))^{C_2}\right).
\end{align*}
Hence, partial integration yields
\begin{align*}
S(I)=\phi^\ast(q)\mathcal{A}(\vartheta(n),n,q)&\int\limits_{I} g(t)\mathrm{d}t+E(t_2)g(t_2)-E(t_1)g(t_1)-\int\limits_{I} E(t) g'(t)\mathrm{d}t\\
    =\phi^\ast(q)\mathcal{A}(\vartheta(n),n,q)&\int\limits_{I} g(t)\mathrm{d}t \\
      &+O\left(\Big(\sup_{t\in I}|E(t)|\Big)\Big(|g(t_1)|+|g(t_2)|+\int\limits_{I} |g'(t)|\mathrm{d}t\Big)\right).
\end{align*}
Splitting $I$ into $R(g)+1$ intervals where the sign of $g'$ does not change, we complete this lemma.
\end{proof}

\begin{remark}\label{R1}
In particular, we apply Lemma \ref{L11} with $q=1$, $m=1$ to deal with the linear congruence condition, where
$\mathcal{A}(\vartheta(n),n):=\mathcal{A}(\vartheta(n),n,1)$. While we use Lemma \ref{L11} with $m=2$ and some idea of \cite[Lemma 6.3]{DEFR1} to handle the quadratic congruence condition.
\end{remark}
\subsection{A general preliminary result for step (3)}
$\newline$In order to achieve the above step (3), we need to compute the Peyre constant. It is necessary to determine the value of $\alpha(\widetilde{S})$
which is a factor of it. Thus, we will give a preliminary lemma on $\alpha(\widetilde{S})$ whose proof follows the \cite[Lemma 8.1]{DEFR1}. It may be viewed as the version on $\mathbb{Q}$ of \cite[Lemma 8.1]{DEFR1}. Before proving, we have to quote some notations from \cite[Section 8]{DEFR1}.

Suppose that the negative curves on $\widetilde{S}$ are $E_1, \dots, E_{r+1+s}$, for some $s \geq 0$, where $E_1, \dots, E_{r+1}$ are a basis of
Pic$(\widetilde{S})$. Clearly, $r=9-d$. Expressing $-K_{\widetilde{S}}$ and $E_{r+2},\dots, E_{r+1+s}$ in terms of this basis, we have
\begin{equation*}
 [-K_{\widetilde{S}}]= \sum_{j=1}^{r+1} c_j[E_j]
\end{equation*}
and, for $i=1, \dots, s$,
\begin{equation*}
  [E_{r+1+i}] = \sum_{j=1}^{r+1} b_{i,j}[E_j]
\end{equation*}
for some $b_{i,j}, c_j \in \mathbb{Z}$.

\begin{lemma}\label{L12}
Assume that $c_{r+1} > 0$. For $j = 1, \dots, r$ and $i=1, \dots, s$, define
  \begin{equation*}
    a_{0,j}:= c_j, \quad a_{i,j}:=b_{i,r+1}c_j-b_{i,j}c_{r+1},\quad
    A_0:=1,\quad A_i:=b_{i,r+1}.
  \end{equation*}
  Then
  \begin{equation*}
    \alpha(\widetilde{S})(\log B)^r = \frac{1}{c_{r+1}} \int_{R_1(B)}
    \frac{1}{\eta_1\cdots \eta_r} \mathrm{d}\eta_1\cdots\mathrm{d}\eta_r
  \end{equation*}
  with a domain of integration
  \begin{align*}
    R_1(B)
    :=\left\{(\eta_1,\ldots,\eta_r) \in \mathbb{R}^r_{\geq 1} \,\Bigg|\,
      \prod_{j=1}^r\eta_j^{a_{i,j}} \le B^{A_i}\; \text{for all}\; i \in \{0, \dots, s\}
\right\}.
  \end{align*}
\end{lemma}
\begin{proof}
From the proof of \cite[Lemma 8.1]{DEFR1}, we have
\begin{align*}
 \alpha(\widetilde{S}) = \frac{1}{c_{r+1}} \int_{(t_1, \dots, t_r) \in R_0} \mathrm{d} t_1\cdots \mathrm{d} t_r,
\end{align*}
where
\begin{equation*}
    R_0 := \left\{(t_1, \dots, t_r) \in \mathbb{R}^r_{\geq0}  \,\Bigg|\,
      \text{$\sum_{j=1}^r a_{i,j}t_j \le A_i$ for all $i \in \{0, \dots,
        s\}$}\right\}.
  \end{equation*}
Now the change of coordinates $\eta_i=B^{t_i}$ for $i \in \{1,\dots, r\}$ transforms $R_0$ into $R_1(B)$ and completes this lemma.
\end{proof}

\section{Proof for the quintic del Pezzo surface $S_1$ of type $\mathbf{A}_2$}
\subsection{Passage to a universal torsor}
$\newline$ Let $S_1$, $\widetilde{S}_1$, $U_1$ and $N_{U_1,H}(B)$ be defined in Section 1. The purpose of this section is to establish a bijection between the rational points on the open subset $U_1$ and the integral points on the universal torsor above $\widetilde{S}_1$, which are subject to a number of coprimality conditions. Here, we will follow the strategy in \cite{DETS1} and refer to \cite{DE1}.

Along the way we will introduce new variables $\eta_1,\ldots,\eta_8$. We use the notation
\begin{align*}
\boldsymbol{\eta}:=(\eta_1,\ldots,\eta_5),\qquad\boldsymbol{\eta}':=(\eta_1,\ldots,\eta_6),\qquad\boldsymbol{\eta}'':=(\eta_1,\ldots,\eta_8)
\end{align*}
and, for $(k_1,\ldots,k_5)\in\mathbb{Q}^5$,
\begin{align*}
\boldsymbol{\eta}^{(k_1,k_2,k_3,k_4,k_5)}:=\eta_1^{k_1}\eta_2^{k_2}\eta_3^{k_3}\eta_4^{k_4}\eta_5^{k_5}.
\end{align*}
For $i=1,\ldots,8$, let
\begin{equation*}
  (\mathbb{Z}_i, J_i, J_i')=
  \begin{cases}
    (\mathbb{Z}_{>0}, \mathbb{R}_{\ge 1}, \mathbb{R}), & i=1,\\
   (\mathbb{Z}_{>0}, \mathbb{R}_{\ge 1}, \mathbb{R}_{\ge 1}), & i \in \{2, \dots, 5\},\\
     (\mathbb{Z}_{\neq 0}, \mathbb{R}_{\le -1} \cup \mathbb{R}_{\ge 1}, \mathbb{R}), & i=6,\\
    (\mathbb{Z}, \mathbb{R},\mathbb{R}), & i \in \{7, 8\}.
  \end{cases}
\end{equation*}
In this section, we estimate summations over $\eta_i\in\mathbb{Z}_i$ by integrations over $\eta_i\in J_i$, which we enlarge to $\eta_i\in J_i'$ in Lemma \ref{L5}.

\begin{figure}[ht]
  \centering
  \[\xymatrix@R=0.05in @C=0.05in{ E_7 \ar@{-}[rrr] \ar@1{-}[dr] \ar@1{-}[dd] &&& \li{3} \ar@{-}[dr]\\
      & \li{6} \ar@{-}[r] & \li{5} \ar@{-}[r] & \ex{2} \ar@{-}[r] & \ex{1}\\
      E_8 \ar@{-}[rrr] \ar@{-}[ur] & & & \li{4} \ar@{-}[ur]}\]
  \caption{Configuration of curves on $\widetilde{S}_1$.}
  \label{F1}
\end{figure}
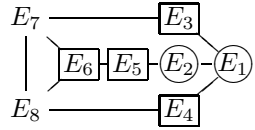

In order to derive the bijection alluded to above, we must begin by collecting
together some useful information about the geometric structure of $S_1$ from \cite[Section 3.3]{DE1}. The singularity $\mathbf{A}_2$ in $(0 : 0 : 1 : 0 : 0 : 0)$ gives exceptional divisors $E_1, E_2$. Let $E_3, E_4, E_5, E_6$ resp. $E_7, E_8$ on $\widetilde{S}_1$
be the strict transforms under $\pi:\widetilde{S}_1\to S_1$ of the four lines $E_3''=\{x_0=x_1=x_3=x_4=0\}$, $E_4''=\{x_0=x_1=x_4=x_5=0\}$, $E_5''=\{x_0=x_3=x_5=x_1+x_4=0\}$, $E_6''=\{x_2=x_3=x_5=x_1+x_4=0\}$ resp. the curves $E_7''=\{x_1=x_2=x_3=x_0x_5+x_4^2=0\}$, $E_8''=\{x_2=x_4=x_5=x_0x_3+x_1^2=0\}$ on $S_1$. The configuration of these divisors $E_1,\ldots,E_8$ on $\widetilde{S}_1$ is described by Figure \ref{F1} with circles and squares marking
classes of $(-2)$- and $(-1)$-curves.

From \cite[Section 3.3]{DE1}, the Cox ring of $\widetilde{S}_1$ is given by
\begin{equation*}
\text{Cox}(\widetilde{S}_1)=\text{Spec}\big(\mathbb{Q}[\boldsymbol{\eta}'']/(T(\boldsymbol{\eta}''))\big),
\end{equation*}
where \begin{equation*}
T(\boldsymbol{\eta}'')=\eta_2\eta_5^2\eta_6+\eta_3\eta_7+\eta_4\eta_8.
\end{equation*}
The universal torsor $\mathcal{T}_1(B)$ is an open subset of the affine hypersurface
\begin{equation}\label{U1}
\eta_2\eta_5^2\eta_6+\eta_3\eta_7+\eta_4\eta_8=0,
\end{equation}
together with a map $\Psi:\mathcal{T}_1(B)\to S_1$, given by
\begin{equation}\label{X1}
\begin{aligned}
(&\Psi^*(x_0),\Psi^*(x_1),\Psi^*(x_2),\Psi^*(x_3),\Psi^*(x_4),\Psi^*(x_5))\\
&=(\boldsymbol{\eta}^{(3,2,2,2,1)},
  \boldsymbol{\eta}^{(2,1,2,1,0)}\eta_7, \eta_6\eta_7\eta_8,\boldsymbol{\eta}^{(1,1,1,0,1)}\eta_6\eta_7, \boldsymbol{\eta}^{(2,1,1,2,0)}\eta_8, \boldsymbol{\eta}^{(1,1,0,1,1)}\eta_6\eta_8).
\end{aligned}
\end{equation}

Next, we will show that the problem of counting rational points of bounded height on the surface $S_1$ can translate into a counting problem for certain
integral points on the universal torsor $\mathcal{T}_1(B)$, subject to coprimality and height inequalities.

\begin{lemma}\label{L1}
 We have
$$N_{U_1,H}(B) = \#\mathcal{T}_1(B),$$
where $\mathcal{T}_1(B)$ is the set of $\boldsymbol{\eta}'' \in \mathbb{Z}_1\times\cdots\times\mathbb{Z}_8$
such that (\ref{U1}) holds, with

\begin{equation}\label{H1}
\max\left\{\begin{aligned}
&|\boldsymbol{\eta}^{(3,2,2,2,1)}|,|\boldsymbol{\eta}^{(2,1,2,1,0)}\eta_7|, |\eta_6\eta_7\eta_8|,\\
&|\boldsymbol{\eta}^{(1,1,1,0,1)}\eta_6\eta_7|,|\boldsymbol{\eta}^{(2,1,1,2,0)}\eta_8|, |\boldsymbol{\eta}^{(1,1,0,1,1)}\eta_6\eta_8|
\end{aligned}\right\} \le B
\end{equation}
\begin{align}\label{C1}
\boldsymbol{\eta}''\; \text{fulfills coprimality conditions as in Figure \ref{F1}}.
\end{align}
\end{lemma}

Coprimality conditions are described by Figure \ref{F1} in the following sense. Let $E_i$ corresponds to $\eta_i$ with $1\leq i \leq 8$. Then, two coordinates are required to be coprime if and only if the corresponding vertices in Figure \ref{F1} are not connected by an edge.
\begin{proof}
 It is based on the construction of the minimal desingularization $\pi: \widetilde{S}_1 \to S_1$ by the following blow-ups
 $\rho: \widetilde{S}_1 \to \mathbb{P}^2$ in four points. Starting with the curves
\begin{equation*}
 E^{(0)}_1:= \{y_0=0\},\;E^{(0)}_7:=\{y_1=0\},\; E^{(0)}_8:=\{y_2 = 0\},\; E^{(0)}_6:=\{-y_1-y_2 = 0\}
\end{equation*}
in $\mathbb{P}^2$:
  \begin{enumerate}
  \item blow up $E^{(0)}_1\cap E^{(0)}_6$, giving $E_2^{(1)}$,
  \item blow up $E^{(1)}_2\cap E^{(1)}_6$, giving $E^{(2)}_5$,
  \item blow up $E^{(2)}_1\cap E^{(2)}_7$, giving $E^{(3)}_3$,
  \item blow up $E^{(3)}_1\cap E^{(3)}_8$, giving $E^{(4)}_4$.
  \end{enumerate}

We may choose $\rho(E_1), \rho(E_7), \rho(E_8)$ as the coordinate lines in $\mathbb{P}^2 = \{(\eta_1: \eta_7: \eta_8)\}$. Then $\rho(E_6)$
is the line $\eta_7+\eta_8=0$. The morphisms $\rho, \pi$ and the projection
\begin{equation*}
  \begin{array}[h]{cccc}
    \phi: & S_1 & \dashrightarrow & \mathbb{P}^2,\\
    & \mathbf{x} & \mapsto & (x_0:x_1:x_4)
  \end{array}
\end{equation*}
form a following commutative diagram of rational maps between $S_1, \widetilde{S}_1$ and $\mathbb{P}^2$.
\begin{equation*}
  \xymatrix{\widetilde{S}_1 \ar@{->}^\rho[dr]\ar@{->}^\pi[d] & \\
    S_1\ar@{-->}^\phi[r] & \mathbb{P}^2}
\end{equation*}
The inverse map of $\phi$ is
\begin{equation}\label{P1}
  \begin{array}[h]{cccc}
    \phi': & \mathbb{P}^2&\dashrightarrow &S_1,\\
    &(\eta_1: \eta_7: \eta_8)&\mapsto&(\eta_1^3:\eta_1^2\eta_7:\eta_6\eta_7\eta_8:\eta_1\eta_6\eta_7:\eta_1^2\eta_8:\eta_1\eta_6\eta_8)
  \end{array}
\end{equation}
where $\eta_6=-\eta_7-\eta_8$. The maps $\phi, \phi^\prime$ give a bijection between the complement $U_1$ of the lines on $S_1$ and $\{(\eta_1: \eta_7: \eta_8) \in \mathbb{P}^2\mid \gcd(\eta_1,\eta_7,\eta_8)=1\}$, and furthermore, induce a bijection between $U_1(\mathbb{Q})$ and the integral points
\begin{equation*}
  \{(\eta_1, \eta_6, \eta_7, \eta_8) \in \mathbb{Z}_{>0}\times\mathbb{Z}_{\neq 0}\times \mathbb{Z}^2 \mid
  \gcd(\eta_1, \eta_6, \eta_8)=1,\ \eta_6+\eta_7+\eta_8=0\}.
\end{equation*}
Since the above four blow-ups (1)-(4), we have the following steps:

1) Let $\eta_2:=\gcd(\eta_1, \eta_6)\in\mathbb{Z}_{>0}$. Then
\begin{equation*}
\eta_1 =\eta_2\eta_1',\qquad \eta_6= \eta_2\eta_6'\qquad\text{with}\, \gcd(\eta_1', \eta_6') = 1.
\end{equation*}
After renaming the variables, we have
\begin{equation*}
T_1 = \eta_2\eta_6+\eta_7+\eta_8=0
\end{equation*}
and
\begin{equation*}
\phi_1': (\eta_1 : \eta_2 : \eta_6 : \eta_7 : \eta_8)\mapsto(\eta_1^3\eta_2^2:\eta_1^2\eta_2\eta_7:\eta_6\eta_7\eta_8:\eta_1\eta_2\eta_6\eta_7:\eta_1^2\eta_2\eta_8:\eta_1\eta_2\eta_6\eta_8).
\end{equation*}

2) Let $\eta_5:=\gcd(\eta_2, \eta_6)\in\mathbb{Z}_{>0}$. Then
\begin{equation*}
\eta_2 =\eta_5\eta_2',\qquad \eta_6= \eta_5\eta_6'\qquad\text{with}\, \gcd(\eta_2', \eta_6') = 1.
\end{equation*}
After renaming the variables, we have
\begin{equation*}
T_2 = \eta_2\eta_5^2\eta_6+\eta_7+\eta_8=0
\end{equation*}
and
\begin{align*}
\phi_2': (\eta_1 : \eta_2 : \eta_5& : \eta_6 : \eta_7 : \eta_8)\mapsto\\
&(\eta_1^3\eta_2^2\eta_5:\eta_1^2\eta_2\eta_7:\eta_6\eta_7\eta_8:\eta_1\eta_2\eta_5\eta_6\eta_7:\eta_1^2\eta_2\eta_8:\eta_1\eta_2\eta_5\eta_6\eta_8).
\end{align*}

3) Let $\eta_3:=\gcd(\eta_1, \eta_7)\in\mathbb{Z}_{>0}$. Then
\begin{equation*}
\eta_1 =\eta_3\eta_1',\qquad \eta_7= \eta_3\eta_7'\qquad\text{with}\,\gcd(\eta_1', \eta_7') = 1.
\end{equation*}
After renaming the variables, we have
\begin{equation*}
T_3 = \eta_2\eta_5^2\eta_6+\eta_3\eta_7+\eta_8=0
\end{equation*}
and
\begin{align*}
\phi_3': (\eta_1 : \eta_2 :& \eta_3 : \eta_5 : \eta_6 : \eta_7: \eta_8)\mapsto\\
&(\eta_1^3\eta_2^2\eta_3^2\eta_5:\eta_1^2\eta_2\eta_3^2\eta_7:\eta_6\eta_7\eta_8:\eta_1\eta_2\eta_3\eta_5\eta_6\eta_7:
\eta_1^2\eta_2\eta_3\eta_8:\eta_1\eta_2\eta_5\eta_6\eta_8)
\end{align*}

4) Let $\eta_4:=\gcd(\eta_1, \eta_8)\in\mathbb{Z}_{>0}$. Then
\begin{equation*}
\eta_1 =\eta_4\eta_1',\qquad \eta_8=\eta_4\eta_8'\qquad\text{with}\,\gcd(\eta_1', \eta_8') = 1.
\end{equation*}
After renaming the variables, we have
\begin{equation*}
T_4=  \eta_2\eta_5^2\eta_6+\eta_3\eta_7+\eta_4\eta_8=0
\end{equation*}
and
\begin{align*}
\phi_4': &(\eta_1 : \eta_2 : \eta_3 : \eta_4 : \eta_5 : \eta_6 : \eta_7: \eta_8)\mapsto\\
&(\eta_1^3\eta_2^2\eta_3^2\eta_4^2\eta_5:\eta_1^2\eta_2\eta_3^2\eta_4\eta_7:\eta_6\eta_7\eta_8:\eta_1\eta_2\eta_3\eta_5\eta_6\eta_7:
\eta_1^2\eta_2\eta_3\eta_4^2\eta_8:\eta_1\eta_2\eta_4\eta_5\eta_6\eta_8)
\end{align*}

Obviously, at each stage the coprimality conditions correspond to intersection properties of the respective divisors. The final result is summarized in
Figure \ref{F1}, which encodes data from (\ref{C1}). Note that $\phi_4'$ is  $\Psi$ from (\ref{X1}). Therefore, $H(\phi_4'(\boldsymbol{\eta}'')) \le B$ is equivalent to (\ref{H1}).
\end{proof}

\subsection{The first summation over $\eta_7$}
$\newline$ From Lemma \ref{L1}, our counting problem has the special form of \cite[Table 1]{DE2}. Our Table \ref{table1} provides a dictionary between the notation of \cite[Section 2]{DE2} and the present situation on $S_1$. We suggest that the readers can refer to \cite[(2.1)]{DE2} and \cite[Definition 2.2]{DE2} to understand the notation of our Table \ref{table1}.

\begin{table}[!ht]
\centering
\caption{Dictionary for applying \cite[Proposition 2.4]{DE2}}
\begin{tabular}{|c|c|c|c|} \hline
$(r,s,t)$ & $(2,1,1)$& $ \delta$ & $\eta_1$\\ \hline
($\alpha_0$;$\alpha_1,\ldots,\alpha_r$) & $(\eta_6;\eta_2,\eta_5)$ &
$(a_0;a_1,\ldots,a_r)$ & $(1;1,2)$ \\ \hline
($\beta_0$;$\beta_1,\ldots,\beta_s$) & $(\eta_7;\eta_3)$ &
$(b_0;b_1,\ldots,b_s)$ & $(1;1)$ \\ \hline
($\gamma_0$;$\gamma_1,\ldots,\gamma_t$) & $(\eta_8;\eta_4)$ &
$(c_1,\ldots,c_t)$ & $(1)$ \\ \hline
$\Pi(\boldsymbol{\alpha})$ & $\eta_2\eta_5^2$ & $\Pi'(\delta,\boldsymbol{\alpha})$ & $\eta_1\eta_2$\\ \hline
$\Pi(\boldsymbol{\beta})$ & $\eta_3$ & $\Pi'(\delta,\boldsymbol{\beta})$ & $\eta_1$\\ \hline
$\Pi(\boldsymbol{\gamma})$ & $\eta_4$ & $\Pi'(\delta,\boldsymbol{\gamma})$ & $\eta_1$\\
\hline
\end{tabular}
\label{table1}
\end{table}

Using (\ref{U1}) to eliminate $\eta_8$, from the height condition (\ref{H1}), we define $\mathcal{R}(B)$ to be the set of all $(\boldsymbol{\eta}',\eta_7) \in
J_1\times\cdots \times J_7$ and
\begin{align}
\left|\boldsymbol{\eta}^{(3,2,2,2,1)}\right|&\leq B, \label{H11}\\
\left|\boldsymbol{\eta}^{(2,1,2,1,0)}\eta_7\right|&\leq B, \label{H12}\\
\left|\frac{\eta_2\eta_5^2\eta_6^2\eta_7+\eta_3\eta_6\eta_7^2}{\eta_4}\right|&\leq B, \label{H13}\\
\left|\boldsymbol{\eta}^{(1,1,1,0,1)}\eta_6\eta_7\right|&\leq B, \label{H14}\\
\left|\boldsymbol{\eta}^{(2,2,1,1,2)}\eta_6+\boldsymbol{\eta}^{(2,1,2,1,0)}\eta_7\right|&\leq B, \label{H15}\\
\left|\boldsymbol{\eta}^{(1,2,0,0,3)}\eta_6^2+\boldsymbol{\eta}^{(1,1,1,0,1)}\eta_6\eta_7\right|&\leq B. \label{H16}
\end{align}
From (\ref{H12}), (\ref{H14}), (\ref{H15}) and (\ref{H16}), we can imply
\begin{align}
\left|\boldsymbol{\eta}^{(2,2,1,1,2)}\eta_6\right|&\leq 2B, \label{H155}\\
\left|\boldsymbol{\eta}^{(1,2,0,0,3)}\eta_6^2\right|&\leq 2B. \label{H166}
\end{align}

Let $\vartheta_0(\boldsymbol{\eta}'):=\prod\limits_{p}\vartheta_{0,p}\left(I_p(\boldsymbol{\eta}')\right)$
with $I_p(\boldsymbol{\eta}')=\left\{i\in\{1,\ldots,6\}:p\mid \eta_i\right\}$ and
\begin{align*}
\vartheta_{0,p}(I)=
\begin{cases}
     1,  &I=\emptyset, \{1\}, \{2\}, \{3\}, \{4\}, \{5\}, \{6\},
    \{1,3\}, \{1,2\}, \{2,5\}, \{5,6\}, \{1,4\},\\
    0, & \text{all other} \;I \subset \{1, \dots, 6\}.
\end{cases}
\end{align*}
Then
\begin{align}\label{Special1}
\vartheta_0(\boldsymbol{\eta}')=1
\end{align}
if and only if $\boldsymbol{\eta}'$ fulfills coprimality conditions from (\ref{C1}).

\begin{lemma}\label{L2}
We have
$$N_{U_1,H}(B) = \sum_{\substack{\boldsymbol{\eta}' \in \mathbb{Z}_1\times\cdots\times\mathbb{Z}_6}}
    \vartheta_1(\boldsymbol{\eta}')V_1(\boldsymbol{\eta}';B) + O\left(B(\log B)^2\right),$$
where
\begin{align}\label{E1}
V_1(\boldsymbol{\eta}';B)=\frac{1}{\eta_4}\int_{(\boldsymbol{\eta}', \eta_7)\in\mathcal{R}(B)}\mathrm{d}\eta_7
\end{align}
and
\begin{align}\label{E2}
\vartheta_1(\boldsymbol{\eta}')=\prod\limits_{p}\vartheta_{1,p}\left(I_p(\boldsymbol{\eta}')\right)
\end{align}
with $I_p(\boldsymbol{\eta}')=\left\{i\in\{1,\ldots,6\}:p\mid \eta_i\right\}$ and
\begin{align*}
\vartheta_{1,p}(I)=
\begin{cases}
     1,  &I=\emptyset, \{3\}, \{4\}, \{6\},\\
    1-\frac{1}{p},&I=\{2\}, \{5\}, \{1,3\}, \{1,2\}, \{2,5\}, \{5,6\}, \{1,4\},\\
    1-\frac{2}{p}, &I =\{1\},\\
    0, & \text{all other} \;I \subset \{1, \dots, 6\}.
\end{cases}
\end{align*}
\end{lemma}
\begin{proof}
From Table \ref{table1} and \cite[Proposition 2.4]{DE2}, we have
\begin{align*}
N_{U_1,H}(B) =\sum_{\substack{\boldsymbol{\eta}' \in\mathbb{Z}_1\times\cdots\times\mathbb{Z}_6}}
     \Big(\vartheta_1(\boldsymbol{\eta}')V_1(\boldsymbol{\eta}';B) +R_1(\boldsymbol{\eta}';B)\Big),
\end{align*}
where $V_1(\boldsymbol{\eta}';B)$ is defined as (\ref{E1}),
\begin{align}\label{Linear}
\vartheta_1(\boldsymbol{\eta}')&=\sum\limits_{\substack{k\mid \eta_1\\ \gcd(k,\eta_2\eta_3\eta_5\eta_6)=1}}
 \frac{\mu(k)\phi^\ast(\eta_1\eta_2\eta_5)}{k\phi^\ast(\gcd(\eta_1,k\eta_4))}
    \sum\limits_{\substack{\eta_2\eta_5^2\eta_6\equiv-\varrho\eta_3\pmod{k\eta_4}\\
    1 \le \varrho \le k\eta_4\\
    \gcd(\varrho,k\eta_4)=1}}1
\end{align}
and\begin{align*}
R_1(\boldsymbol{\eta}';B)\ll 2^{\omega(\eta_1)+\omega(\eta_1\eta_2\eta_5)}.
\end{align*}

In order to decide the main term, it suffices to compute $\vartheta_1(\boldsymbol{\eta}')$. Similar to \cite[Lemma 5.4]{DEFR1}, we have
\begin{align*}
\vartheta_1(\boldsymbol{\eta}')&=\sum\limits_{\substack{k\mid \eta_1\\ \gcd(k,\eta_2\eta_3\eta_5\eta_6)=1}}
 \frac{\mu(k)\phi^\ast(\eta_1\eta_2\eta_5)}{k\phi^\ast(\gcd(\eta_1,k\eta_4))}\\
 &=\frac{\phi^\ast(\eta_1\eta_2\eta_5)}{\phi^\ast(\gcd(\eta_1,\eta_4))}
 \sum\limits_{\substack{k\mid \eta_1\\ \gcd(k,\eta_2\eta_3\eta_5\eta_6)=1}}\frac{\mu(k)}{k}\prod\limits_{\substack{p\mid\gcd(k,\eta_1)\\p\nmid \eta_4}}\left(1-\frac{1}{p}\right)^{-1}.
\end{align*}
Further, by (\ref{Special1}), we have
\begin{align*}
\frac{\phi^\ast(\eta_1\eta_2\eta_5)}{\phi^\ast(\gcd(\eta_1,\eta_4))}&=\frac{\prod\limits_{p\mid\eta_1\eta_2\eta_5}\left(1-\frac{1}{p}\right)}
{\prod\limits_{p\mid\gcd(\eta_1,\eta_4)}\left(1-\frac{1}{p}\right)},\\
\sum\limits_{\substack{k\mid \eta_1\\ \gcd(k,\eta_2\eta_3\eta_5\eta_6)=1}}\frac{\mu(k)}{k}&=\prod\limits_{\substack{p\mid\eta_1\\p\nmid \eta_2\eta_3\eta_5\eta_6}}\left(1-\frac{1}{p}\right),\\
\prod\limits_{\substack{p\mid\gcd(k,\eta_1)\\p\nmid \eta_4}}\left(1-\frac{1}{p}\right)^{-1}&=
\prod\limits_{\substack{p\mid\eta_1\\p\nmid \eta_2\eta_3\eta_4\eta_5\eta_6}}\frac{p-2}{p-1}\prod\limits_{\substack{p\mid\gcd(\eta_1,\eta_4)\\p\nmid \eta_2\eta_3\eta_5\eta_6}}\left(1-\frac{1}{p}\right).
\end{align*}
Since (\ref{Special1}) and a straightforward inspection of the above local factors, (\ref{E2}) can be proved.

Next, we will utilize (\ref{H155}) to deal with the error term.
\begin{align*}
\sum_{\substack{\boldsymbol{\eta}'}}R_1(\boldsymbol{\eta}';B)
&\ll \sum_{\substack{\boldsymbol{\eta}'}}2^{\omega(\eta_1)+\omega(\eta_1\eta_2\eta_5)}\\
&\ll B\sum_{\substack{\boldsymbol{\eta}}}\frac{2^{\omega(\eta_1)+\omega(\eta_1\eta_2\eta_5)}}{\boldsymbol{\eta}^{(2,2,1,1,2)}}\\
&\ll B(\log B)^2.
\end{align*}
\end{proof}

\subsection{The second summation over $\eta_6$}
\begin{lemma}\label{L3}
We have
\begin{align*}
N_{U_1,H}(B) = \sum_{\substack{\boldsymbol{\eta} \in \mathbb{Z}_1\times\cdots\times\mathbb{Z}_5}}
   \mathcal{A}(\vartheta_1(\boldsymbol{\eta}'),\eta_6)V_{11}(\boldsymbol{\eta};B) + O\left(B(\log B)^3\log\log B\right),
\end{align*}
where $\mathcal{A}(\vartheta_1(\boldsymbol{\eta}'),\eta_6)$ is defined in Section 2 and
\begin{align}\label{V1}
V_{11}(\boldsymbol{\eta};B) =\frac{1}{\eta_4}\int_{(\boldsymbol{\eta}', \eta_7)\in\mathcal{R}(B)}\mathrm{d}\eta_6\mathrm{d}\eta_7.
\end{align}
\end{lemma}
\begin{proof}
Firstly, we rewrite the result of Lemma \ref{L2} as follow.
\begin{align*}
N_{U_1,H}(B) = \sum_{\substack{\boldsymbol{\eta} \in \mathbb{Z}_1\times\cdots\times\mathbb{Z}_5}}\sum_{\eta_6\in\mathbb{Z}_6}
    \vartheta(\eta_6)g(\eta_6) + O\left(B(\log B)^2\right),
\end{align*}
where\begin{align*}
\vartheta(\eta_6):= \vartheta_1(\boldsymbol{\eta},\eta_6) \qquad\text{and}\qquad g(\eta_6):= V_1(\boldsymbol{\eta},\eta_6;B).
\end{align*}

Since (\ref{E2}) and Definition \ref{D1}, we have $c=1$ and $A_p(0)=1$. Further, by Definition \ref{D2}, we yield that
\begin{align*}
\vartheta(\eta_6)\in \Theta_2(b,1,1,1)\qquad\text{and}\qquad b=\prod\limits_{p\mid \eta_1\eta_2\eta_3\eta_4}p.
\end{align*}
Applying \cite[Lemma 5.1(4)]{DE2} to (\ref{H13}), we imply that
\begin{align}\label{E4}
g(\eta_6)\ll\left(\frac{B}{\eta_3\eta_4|\eta_6|}\right)^\frac{1}{2}
=\frac{B}{\boldsymbol{\eta}^{(1,1,1,1,1)}|\eta_6|}\left(\frac{B}{\boldsymbol{\eta}^{(3,2,2,2,1)}}\right)^{-1/4}
\left(\frac{B}{\boldsymbol{\eta}^{(1,2,0,0,3)}\eta_6^2}\right)^{-1/4}.
\end{align}
In particular, by (\ref{H11}) and (\ref{H166}), we deduce that
\begin{align}\label{E5}
g(\eta_6)\ll\frac{B}{\boldsymbol{\eta}^{(1,1,1,1,1)}|\eta_6|}.
\end{align}

Let $t_1 := (\log B)^7$. Since a straightforward application of Lemma \ref{L11} would not yield sufficiently good error terms, we will
use a strategy from \cite[Lemma 4]{DEFR2}. Thus, we split the sum over $\eta_6$ into the two cases $|\eta_6| \leq t_1$ and $|\eta_6| > t_1$.

For the second case, from (\ref{Linear}), we learn that the sum over $\varrho$ is linear and is just 1.
By Lemma \ref{L11}, Remark \ref{R1} and (\ref{E5}), we get that
\begin{align*}
\sum_{\substack{\eta_6\in\mathbb{Z}_6\\|\eta_6| > t_1}}\vartheta(\eta_6)g(\eta_6)=\mathcal{A}(\vartheta(\eta_6),\eta_6)\int\limits_{\substack{\eta_6\in\mathbb{Z}_6\\|\eta_6| > t_1}}g(\eta_6)\mathrm{d}\eta_6
+O\left(\frac{B\log B2^{\omega(\eta_1\eta_2\eta_3\eta_4)}}{\boldsymbol{\eta}^{(1,1,1,1,1)}|\eta_6|}\right),
\end{align*}
where we have used the fact
\begin{align*}
\tau\left(\prod\limits_{p\mid \eta_1\eta_2\eta_3\eta_4}p\right)=2^{\omega(\eta_1\eta_2\eta_3\eta_4)}.
\end{align*}
When summing the error term over the remaining variables with $|\eta_6| > t_1$, we have
\begin{align*}
\ll \frac{B\log B}{|\eta_6|}\sum_{\substack{\boldsymbol{\eta}}}\frac{2^{\omega(\eta_1\eta_2\eta_3\eta_4)}}{\boldsymbol{\eta}^{(1,1,1,1,1)}}
\ll t_1^{-1}B(\log B)^{10}\ll B(\log B)^3.
\end{align*}

Now, we start with the first case. Since $|\eta_6|\leq t_1$ and $0\leq\vartheta(\eta_6)\leq 1$, by (\ref{E4}), we yield an upper bound
\begin{equation}\label{E6}
\begin{aligned}
\sum_{\substack{\boldsymbol{\eta} \in \mathbb{Z}_1\times\cdots\times\mathbb{Z}_5}}
\sum_{\substack{\eta_6\in\mathbb{Z}_6\\|\eta_6| \leq t_1}}&\vartheta(\eta_6)g(\eta_6)\\
&\ll \sum_{\substack{\boldsymbol{\eta}'\\ |\eta_6| \leq t_1}}\frac{B}{\boldsymbol{\eta}^{(1,1,1,1,1)}|\eta_6|}\left(\frac{B}{\boldsymbol{\eta}^{(3,2,2,2,1)}}\right)^{-1/4}
\left(\frac{B}{\boldsymbol{\eta}^{(1,2,0,0,3)}\eta_6^2}\right)^{-1/4}\\
&\ll \sum_{\substack{\eta_1,\eta_2,\eta_4,\eta_5,\eta_6\\ |\eta_6| \leq t_1}}\frac{B}{\eta_1\eta_2\eta_4\eta_5|\eta_6|}\left(\frac{B}{\boldsymbol{\eta}^{(1,2,0,0,3)}\eta_6^2}\right)^{-1/4}\\
&\ll  \sum_{\substack{\eta_1,\eta_2,\eta_4,\eta_6\\ |\eta_6| \leq t_1}}\frac{B}{\eta_1\eta_2\eta_4|\eta_6|}\ll B(\log B)^3\log t_1\\
&\ll B(\log B)^3\log\log B,
\end{aligned}
\end{equation}
where we have used the conditions (\ref{H11}) and (\ref{H166}).

In order to prove this lemma, it suffices to obtain the following upper bound.
\begin{align*}
\sum_{\substack{\boldsymbol{\eta} \in \mathbb{Z}_1\times\cdots\times\mathbb{Z}_5}}&\mathcal{A}(\vartheta(\eta_6),\eta_6)\int\limits_{\substack{\eta_6\in\mathbb{Z}_6\\|\eta_6| \leq t_1}}g(\eta_6)\mathrm{d}\eta_6\\
&\ll\sum_{\substack{\boldsymbol{\eta}}}\int\limits_{|\eta_6| \leq t_1}\frac{B}{\boldsymbol{\eta}^{(1,1,1,1,1)}|\eta_6|}
\left(\frac{B}{\boldsymbol{\eta}^{(3,2,2,2,1)}}\right)^{-1/4}\left(\frac{B}{\boldsymbol{\eta}^{(1,2,0,0,3)}\eta_6^2}\right)^{-1/4}\mathrm{d}\eta_6\\
&\ll\sum_{\substack{\eta_1,\eta_2,\eta_4}}\frac{B}{\eta_1\eta_2\eta_4}\int\limits_{|\eta_6| \leq t_1}|\eta_6|^{-1}\mathrm{d}\eta_6
\ll B(\log B)^3\log\log B,
\end{align*}
where we have used (\ref{E6}), $0\leq \mathcal{A}(\vartheta(\eta_6),\eta_6)\leq 1$.

We complete this lemma, since
\begin{align*}
\mathcal{A}(\vartheta(\eta_6),\eta_6)=\mathcal{A}(\vartheta_1(\boldsymbol{\eta}'),\eta_6)\qquad\text{and}\qquad
\int\limits_{\substack{\eta_6\in\mathbb{Z}_6}} g(\eta_6)\mathrm{d}\eta_6=V_{11}(\boldsymbol{\eta};B).
\end{align*}

\end{proof}

\subsection{The remaining summations}
\begin{lemma}\label{L4}
We have
\begin{align*}
N_{U_1,H}(B)=\left(\prod\limits_{p}\omega_p\right)V_0(B)+O\left(B(\log B)^3\log\log B\right),
\end{align*}
where
\begin{align*}
\omega_p:=\left(1-\frac{1}{p}\right)^5\left(1+\frac{5}{p}+\frac{1}{p^2}\right)
\end{align*} and
\begin{align*}
V_0(B): =\int_{(\boldsymbol{\eta}', \eta_7)\in\mathcal{R}(B)}\frac{1}{\eta_4}\mathrm{d}\eta_1\cdots\mathrm{d}\eta_7.
\end{align*}
\end{lemma}
\begin{proof}
Firstly, we rewrite the result of Lemma \ref{L3} as follow.
\begin{align}\label{V2}
N_{U_1,H}(B) =\sum_{\substack{\boldsymbol{\eta} \in \mathbb{Z}_1\times\cdots\times\mathbb{Z}_5}}
   \vartheta_5(\boldsymbol{\eta})V_5(\boldsymbol{\eta};B)+O\left(B(\log B)^3\log\log B\right),
\end{align}
where\begin{align*}
\vartheta_5(\boldsymbol{\eta}):= \mathcal{A}(\vartheta_1(\boldsymbol{\eta}'),\eta_6)\qquad\text{and}\qquad
V_5(\boldsymbol{\eta};B):= V_{11}(\boldsymbol{\eta};B).
\end{align*}

Since (\ref{E2}), we have $\vartheta_1(\boldsymbol{\eta}')\in\Theta'_{4,6}(2)$ \cite[Definition 7.8]{DE2}. By \cite[Fig.2]{DE2}, we get $\vartheta_1(\boldsymbol{\eta}')\in\Theta_{2,6}(C)$ \cite[Definition 4.2]{DE2} for some $C\in\mathbb{Z}_{\geq0}$. Furthermore, we have
$\vartheta_5(\boldsymbol{\eta})=\mathcal{A}(\vartheta_1(\boldsymbol{\eta}'),\eta_6)\in\Theta_{2,5}(C)$. Applying \cite[Lemma 5.1(5)]{DE2} to (\ref{H13}), we have
\begin{align*}
V_5(\boldsymbol{\eta};B)\ll\frac{B^{2/3}}{(\eta_2\eta_3\eta_4\eta_5^2)^{1/3}}
=\frac{B}{\boldsymbol{\eta}^{(1,1,1,1,1)}}\left(\frac{B}{\boldsymbol{\eta}^{(3,2,2,2,1)}}\right)^{-1/3}.
\end{align*}
Now, the conditions of \cite[Proposition 4.3]{DE2} are satisfied, and we have
\begin{equation}\label{A}
\begin{aligned}
& \sum_{\substack{\eta_1,\ldots,\eta_{r+s}}}\vartheta_{r+s}(\eta_1,\ldots,\eta_{r+s})V_{r+s}(\eta_1,\ldots,\eta_{r+s};B)\\
=&\mathcal{A}(\vartheta_{r+s}(\eta_1,\ldots,\eta_{r+s}),\eta_1,\ldots,\eta_{r+s}))
   \int\limits_{\eta_1,\ldots,\eta_{r+s}}V_{r+s}(\eta_1,\ldots,\eta_{r+s};B)\mathrm{d}\eta_1\cdots\mathrm{d}\eta_{r+s}\\
   &+O\left(B(\log B)^{r-1}(\log\log B)^{\max\{1,s\}}\right).
\end{aligned}
\end{equation}
Applying this to the main term of (\ref{V2}) with $(r,s)=(4,1)$, we imply that
\begin{align*}
\sum_{\substack{\boldsymbol{\eta} \in \mathbb{Z}_1\times\cdots\times\mathbb{Z}_5}}
   \vartheta_5(\boldsymbol{\eta})V_5(\boldsymbol{\eta};B)=\mathcal{A}( \vartheta_5(\boldsymbol{\eta}),\boldsymbol{\eta}) \int\limits_{\boldsymbol{\eta}}V_5(\boldsymbol{\eta};B)\mathrm{d}\boldsymbol{\eta}+O\left(B(\log B)^3\log\log B\right).
\end{align*}
By (\ref{V1}), we have
\begin{align*}
\int\limits_{\boldsymbol{\eta}}V_5(\boldsymbol{\eta};B)\mathrm{d}\boldsymbol{\eta}=V_0(B).
\end{align*}

Since $\mathcal{A}( \vartheta_5(\boldsymbol{\eta}),\boldsymbol{\eta})=\mathcal{A}(\vartheta_1(\boldsymbol{\eta}'),\boldsymbol{\eta}')$ and
$\vartheta_1(\boldsymbol{\eta}')\in\Theta'_{4,6}(2)$, we get that the conditions of \cite[Corollary 7.10]{DE2} are satisfied.
From \cite[Corollary 7.10]{DE2} and (\ref{E2}), we have
\begin{align*}
\mathcal{A}( \vartheta_5(\boldsymbol{\eta}),\boldsymbol{\eta})=\prod\limits_{p}\omega_p,
\end{align*}
where
\begin{align*}
\omega_p=&\left(1-\frac{1}{p}\right)^6+3\left(1-\frac{1}{p}\right)^5\left(\frac{1}{p}\right)
+2\left(1-\frac{1}{p}\right)^5\left(\frac{1}{p}\right)\left(1-\frac{1}{p}\right)\\
&+5\left(1-\frac{1}{p}\right)^4\left(\frac{1}{p}\right)^2\left(1-\frac{1}{p}\right)
+\left(1-\frac{1}{p}\right)^5\left(\frac{1}{p}\right)\left(1-\frac{2}{p}\right)\\
=&\left(1-\frac{1}{p}\right)^5\left(1+\frac{5}{p}+\frac{1}{p^2}\right).
\end{align*}
\end{proof}

\subsection{The final proof}
$\newline$For a split singular del Pezzo surface $S_1$, by  \cite[Theorem 4]{DE4} and \cite[Theorem 1.3]{DEJOTE1}, it is easy to compute the constant $\alpha(\widetilde{S}_1)$. We have
\begin{align}\label{W1}
\alpha(\widetilde{S}_1)=\frac{1}{144}\cdot\frac{1}{\#W(\mathbf{A}_2)}=\frac{1}{864},
\end{align}
where $W(\mathbf{A}_n)$ stands for the Weyl group associated to the Dynkin diagram of the singularity $\mathbf{A}_n$. Here, we have used $\#W(\mathbf{A}_n)=(n + 1)!$. Let
\begin{equation*}
\omega_\infty(\widetilde{S}_1):=3\int_{\max\{|z_0^3|,
|z_0^2z_1|, |z_1z_2(z_1+z_2)|, |z_0z_1(z_1+z_2)|, |z_0^2z_2|, |z_0z_2(z_1+z_2)|\}\leq 1}\mathrm{d}z_0 \mathrm{d}z_1 \mathrm{d}z_2.
\end{equation*}
For $(\eta_2,\eta_3,\eta_4,\eta_5)\in\mathbb{R}^4_{\geq 1}$, we introduce the conditions
\begin{align}
&\eta_2^2\eta_3^2\eta_4^2\eta_5\leq B,\label{H17}\\
&\eta_2^2\eta_3^2\eta_4^2\eta_5\leq B,\qquad \eta_2^2\eta_3^{-1}\eta_4^{-1}\eta_5^4\leq B.\label{H18}
\end{align}

\begin{lemma}\label{L5}
 Let $\alpha(\widetilde{S}_1)$ be as in (\ref{W1}), we have
  \begin{equation*}
    V_0'(B) := \int_{\substack{(\boldsymbol{\eta}', \eta_7)\in \mathcal{R}'(B)\\
    (\ref{H18})}} \frac{1}{\eta_4} \mathrm{d}\eta_1\cdots\mathrm{d}\eta_7,
  \end{equation*}
where\begin{align*}
\mathcal{R}'(B)=\{(\boldsymbol{\eta}', \eta_7)\in J_1'\times\cdots\times J_7'\mid(\ref{H11})-(\ref{H16})\}.
\end{align*}
Then $\frac{1}{864}\omega_\infty(\widetilde{S}_1)B(\log B)^4 = V_0'(B)$.
\end{lemma}
\begin{proof}
From (\ref{P1}), we have
\begin{equation}\label{Y1}
\begin{aligned}
(y_0: y_1: y_2)&\mapsto\\
&\left(y_0^3: y_0^2y_1: -y_1y_2(y_1+y_2): -y_0y_1(y_1+y_2): y_0^2y_2: -y_0y_2(y_1+y_2)\right).
\end{aligned}
\end{equation}
Inserting $y_0 = \eta_1\eta_2\eta_3\eta_4\eta_5$, $y_1 =\eta_3\eta_7$, $y_2= \eta_4\eta_8$, $-(y_1+y_2) =
\eta_2\eta_5^2\eta_6$ into (\ref{Y1}) and cancelling out $\eta_2\eta_3\eta_4\eta_5^2$ gives $\phi_4'$ as
in the proof of Lemma \ref{L1}. Let $\eta_2, \eta_3, \eta_4, \eta_5\in\mathbb{R}_{\neq 0}$ and $\eta_1, \eta_6, \eta_7\in\mathbb{R}$. With $l=B\eta_2\eta_3\eta_4\eta_5^2$, we utilize the coordinate transformation
\begin{align*}
z_0=l^{-1/3}\eta_2\eta_3\eta_4\eta_5\cdot\eta_1,\qquad z_1=l^{-1/3}\eta_3\cdot\eta_7,\qquad z_2=l^{-1/3}(-\eta_2\eta_5^2\cdot\eta_6-\eta_3\cdot\eta_7)
\end{align*}
of Jacobi determinant
\begin{align*}
\frac{\eta_2\eta_3\eta_4\eta_5}{B}\cdot\frac{1}{\eta_4}
\end{align*}
and get
\begin{align}\label{W2}
\omega_\infty(\widetilde{S}_1)=\frac{3\eta_2\eta_3\eta_4\eta_5}{B}\int_{(\boldsymbol{\eta}', \eta_7)\in \mathcal{R}''(B)}\frac{1}{\eta_4}\mathrm{d}\eta_1 \mathrm{d}\eta_6 \mathrm{d}\eta_7,
\end{align}
where
\begin{align*}
\mathcal{R}''(B)=\{(\boldsymbol{\eta}', \eta_7)\in\mathbb{R}\times\mathbb{R}^4_{\neq 0}\times\mathbb{R}^2\mid(\ref{H11})-(\ref{H16})\}.
\end{align*}

An application of Lemma \ref{L12} with exchanged roles of $\eta_1$ and $\eta_5$ gives
\begin{align}\label{a1}
\alpha(\widetilde{S}_1)(\log B)^4=\frac{1}{3}\int\limits_{\substack{(\eta_2,\eta_3,\eta_4,\eta_5)\in\mathbb{R}^4_{\geq 1}\\ (\ref{H18})}}
\frac{\mathrm{d}\eta_2\cdots\mathrm{d}\eta_5}{\eta_2\eta_3\eta_4\eta_5},
\end{align}
since $[-K_{\widetilde{S}_1}] = [3E_1+2E_2+2E_3+2E_4+E_5]$, $[E_6]=[E_1+E_3+E_4-E_5]$. From (\ref{W1}) and (\ref{W2})-(\ref{a1}), we have
\begin{align*}
\frac{1}{864}\omega_\infty(\widetilde{S}_1) B(\log B)^4=V_0'(B).
\end{align*}
\end{proof}

To finish our proof, we compare $V_0(B)$ from Lemma \ref{L4} with $V_0'(B)$ from Lemma \ref{L5}.
\begin{lemma}\label{L6}
 We have
\begin{equation*}
V_0(B) = V_0'(B)+O\left(B(\log B)^3\right).
\end{equation*}
\end{lemma}
\begin{proof}
 Let
\begin{align*}
  \mathcal{D}_0(B) &:= \{(\boldsymbol{\eta}', \eta_7) \in \mathbb{R}^7 \, \text{and}\,(\ref{H11})-(\ref{H16})
  \mid \eta_1,\ldots,\eta_5,|\eta_6| \geq 1\},\\
  \mathcal{D}_1(B) &:= \{(\boldsymbol{\eta}', \eta_7) \in \mathbb{R}^7 \, \text{and}\,(\ref{H11})-(\ref{H16})
   \mid \eta_1,\ldots,\eta_5,|\eta_6| \geq 1, (\ref{H17})\},\\
  \mathcal{D}_2(B) &:= \{(\boldsymbol{\eta}', \eta_7) \in \mathbb{R}^7 \, \text{and}\,(\ref{H11})-(\ref{H16})
   \mid \eta_1,\ldots,\eta_5,|\eta_6| \geq 1, (\ref{H18})\},\\
   \mathcal{D}_3(B) &:= \{(\boldsymbol{\eta}', \eta_7) \in \mathbb{R}^7 \, \text{and}\,(\ref{H11})-(\ref{H16})
   \mid \eta_1,\ldots,\eta_5\geq 1, (\ref{H18})\},\\
  \mathcal{D}_4(B) &:= \{(\boldsymbol{\eta}', \eta_7) \in \mathbb{R}^7 \, \text{and}\,(\ref{H11})-(\ref{H16})
   \mid \eta_2,\ldots,\eta_5\geq 1, (\ref{H18})\}.
\end{align*}
For $i\in\{0,\ldots,4\}$, let
\begin{align*}
V_i(B):=\int_{(\boldsymbol{\eta}', \eta_7)\in\mathcal{D}_i(B)}\frac{\mathrm{d}\eta_1\cdots\mathrm{d}\eta_7}{\eta_4}.
\end{align*}
Then $V_0(B)$ is as in Lemma \ref{L4} and $V_4(B) = V_0'(B)$. It suffices to show that $V_i(B)-V_{i-1}(B) = O(B(\log B)^3)$ for $1 \leq i \leq 4$. This holds for $i = 1$, since, by (\ref{H11}) and $\eta_1 \geq 1$, we have $\mathcal{D}_1(B) =\mathcal{D}_0(B)$.

Moreover, using  \cite[Lemma 5.1(4)]{DE2} and (\ref{H13}) to bound the integral over $\eta_7$, we have
\begin{align*}
V_2(B) - V_1(B)& \ll \int\limits_{\substack{\eta_1,\ldots,\eta_5,|\eta_6|\geq 1 \\ \eta_2^2\eta_3^{-1}\eta_4^{-1}\eta_5^4> B\\ (\ref{H166})}}\frac{B^{1/2}}{(\eta_3\eta_4|\eta_6|)^{1/2}}\mathrm{d} \eta_1 \cdots \mathrm{d} \eta_6 \\
&\ll \int\limits_{\substack{\eta_1,\eta_2,\eta_4,\eta_5,|\eta_6|\geq 1\\ (\ref{H166})}}
\frac{\eta_2\eta_5^2}{\eta_4|\eta_6|^{1/2}}\mathrm{d} \eta_1 \mathrm{d} \eta_2\mathrm{d} \eta_4\mathrm{d} \eta_5 \mathrm{d} \eta_6 \\
&\ll \int\limits_{\substack{\eta_2,\eta_4,\eta_5,|\eta_6|\geq 1}}
\frac{B}{\eta_2\eta_4\eta_5|\eta_6|^{5/2}}\mathrm{d} \eta_2\mathrm{d} \eta_4\mathrm{d} \eta_5 \mathrm{d} \eta_6 \\
&\ll B(\log B)^3.
\end{align*}
Similarly, we have
\begin{align*}
V_3(B) - V_2(B)& \ll \int\limits_{\substack{\eta_1,\ldots,\eta_5\geq 1 \\ |\eta_6|<1,(\ref{H11}),(\ref{H18})}}\frac{B^{1/2}}{(\eta_3\eta_4|\eta_6|)^{1/2}}\mathrm{d} \eta_1 \cdots \mathrm{d} \eta_6 \\
&\ll \int\limits_{\substack{\eta_2,\ldots,\eta_5\geq 1 \\ (\ref{H18})}}\frac{B^{5/6}}{\eta_2^{2/3}\eta_3^{7/6}\eta_4^{7/6}\eta_5^{1/3}}\mathrm{d} \eta_2\cdots \mathrm{d} \eta_5\\
&\ll \int\limits_{\substack{\eta_2,\eta_3,\eta_4\geq 1}}\frac{B}{\eta_2\eta_3\eta_4}\mathrm{d} \eta_2\mathrm{d} \eta_3 \mathrm{d} \eta_4\\
&\ll B(\log B)^3.
\end{align*}

Finally, using \cite[Lemma 5.1(5)]{DE2} and (\ref{H13}) to bound the integral over $\eta_6$, $\eta_7$, we have
\begin{align*}
V_4(B)-V_3(B)& \ll \int\limits_{\substack{\eta_2,\ldots,\eta_5\geq 1 \\ \eta_1<1,(\ref{H17})}}\frac{B^{2/3}}{(\eta_2\eta_3\eta_4\eta_5^2)^{1/3}}\mathrm{d} \eta_1 \cdots \mathrm{d} \eta_5 \\
&\ll \int\limits_{\substack{\eta_2,\eta_3,\eta_4\geq 1}}\frac{B}{\eta_2\eta_3\eta_4}\mathrm{d} \eta_2\mathrm{d} \eta_3 \mathrm{d} \eta_4\\
&\ll B(\log B)^3.
\end{align*}
\end{proof}

Theorem 1.1 follows from Lemma \ref{L4}, Lemma \ref{L5} and Lemma \ref{L6}.

\section{Proof for the quartic del Pezzo surface $S_2$ of type $\mathbf{A}_3+\mathbf{A}_1$}
Let $S_2$, $\widetilde{S}_2$, $U_2$ and $N_{U_2,H}(B)$ be defined in Section 1. The proof of Theorem 1.2 is similar to \cite{DE2} and part of steps (1)-(3) of Section 1 on $S_2$ has been achieved in \cite{DE2}.

In detail, for step (1), in \cite[Lemma 8.3]{DE2}, they followed the strategy in \cite{DETS1} to establish a bijection between the rational points on the open subset $U_2$, and the integral points on the universal torsor above $\widetilde{S}_2$. For step (2), in \cite[Lemma 8.4]{DE2}, they used \cite[Proposition 2.4]{DE2} to deal with the first summation over $\eta_8$. Then, in \cite[Lemma 8.5]{DE2}, they applied \cite[Proposition 4.3]{DE2} to deal with the remaining summations $\eta_1,\ldots,\eta_7$. Finally, in \cite[Section 8.3]{DE2}, they completed step (3) in a straightforward manner.

Clearly, they reduced the second summation over $\eta_7$ to the remaining summations in \cite{DE2}. This causes that the error term is not sufficiently good. Therefore, we will use Lemma \ref{L11} and Remark \ref{R1} to deal with the second summation over $\eta_7$ after the first summation over $\eta_8$. Then, we apply \cite[Proposition 4.3]{DE2} to deal with the remaining summations over  $\eta_1,\ldots,\eta_6$.

\subsection{Passage to a universal torsor}
$\newline$We will introduce new variables $\eta_1,\ldots,\eta_9$ and use the notation
\begin{align*}
\boldsymbol{\eta}:=(\eta_1,\ldots,\eta_6),\qquad\boldsymbol{\eta}':=(\eta_1,\ldots,\eta_7),\qquad\boldsymbol{\eta}'':=(\eta_1,\ldots,\eta_9)
\end{align*}
and, for $(k_1,\ldots,k_6)\in\mathbb{Q}^6$,
\begin{align*}
\boldsymbol{\eta}^{(k_1,k_2,k_3,k_4,k_5,k_6)}:=\eta_1^{k_1}\eta_2^{k_2}\eta_3^{k_3}\eta_4^{k_4}\eta_5^{k_5}\eta_6^{k_6}.
\end{align*}
For $i=1,\ldots,9$, let
\begin{equation*}
  (\mathbb{Z}_i, J_i, J_i')=
  \begin{cases}
    (\mathbb{Z}_{>0}, \mathbb{R}_{\ge 1}, \mathbb{R}_{\ge 1}), & i\in \{1, \dots, 5\},\\
   (\mathbb{Z}_{>0}, \mathbb{R}_{\ge 1}, \mathbb{R}_{\ge 0}), & i=6,\\
     (\mathbb{Z}_{\neq 0}, \mathbb{R}_{\le -1} \cup \mathbb{R}_{\ge 1}, \mathbb{R}), & i=7,\\
    (\mathbb{Z}, \mathbb{R},\mathbb{R}), & i \in \{8, 9\}.
  \end{cases}
\end{equation*}
In this section, we estimate summations over $\eta_i\in\mathbb{Z}_i$ by integrations over $\eta_i\in J_i$, which we enlarge to $\eta_i\in J_i'$ in Section 4.5.

\begin{figure}[ht]
  \centering
  \[\xymatrix@R=0.05in @C=0.05in{ E_9 \ar@{-}[rrrr] \ar@1{-}[dr] \ar@1{-}[dd] &&&&  \li{1} \ar@{-}[dr]\\
      & \ex{7} \ar@{-}[r] & \li{5} \ar@{-}[r] & \ex{6}\ar@{-}[r] & \ex{4}\ar@{-}[r] & \ex{3}\\
      E_8 \ar@{-}[rrrr] \ar@{-}[ur] &&&& \li{2} \ar@{-}[ur]}\]
  \caption{Configuration of curves on $\widetilde{S}_2$ (see\cite[Fig.3]{DE2}).}
  \label{F2}
\end{figure}

\begin{lemma}(\cite[Lemma 8.3]{DE2})
\label{L21}
 We have
$$N_{U_2,H}(B) = \#\mathcal{T}_2(B),$$
where $\mathcal{T}_2(B)$ is the set of $\boldsymbol{\eta}'' \in \mathbb{Z}_1\times\cdots\times\mathbb{Z}_9$
such that
\begin{equation}\label{U21}
\eta_1\eta_9+\eta_2\eta_8+\eta_4\eta_5^3\eta_6^2\eta_7=0,
\end{equation}
holds, with
\begin{equation}\label{H21}
\max\left\{\begin{aligned}
&|\boldsymbol{\eta}^{(0,1,1,1,1,1)}\eta_7\eta_8|,|\boldsymbol{\eta}^{(2,2,3,2,0,1)}|, \\ &|\boldsymbol{\eta}^{(1,1,2,2,2,2)}\eta_7|,|\boldsymbol{\eta}^{(0,0,1,2,4,3)}\eta_7^2|,|\eta_7\eta_8\eta_9|
\end{aligned}\right\} \le B,
\end{equation}
\begin{align*}
\boldsymbol{\eta}''\; \text{fulfills coprimality conditions as in Figure \ref{F2}}.
\end{align*}
\end{lemma}

\subsection{The first summation over $\eta_8$}
$\newline$ Using (\ref{U21}) to eliminate $\eta_9$, from the height condition (\ref{H21}), we define $\mathcal{R}(B)$ to be the set of all $(\boldsymbol{\eta}',\eta_8) \in J_1\times\cdots \times J_8$ and
\begin{align}
\left|\boldsymbol{\eta}^{(0,1,1,1,1,1)}\eta_7\eta_8\right|&\leq B, \label{H211}\\
\left|\boldsymbol{\eta}^{(2,2,3,2,0,1)}\right|&\leq B, \label{H212}\\
\left|\boldsymbol{\eta}^{(1,1,2,2,2,2)}\eta_7\right|&\leq B, \label{H213}\\
\left|\boldsymbol{\eta}^{(0,0,1,2,4,3)}\eta_7^2\right|&\leq B, \label{H214}\\
\left|\frac{\eta_2\eta_7\eta_8^2+\eta_4\eta_5^3\eta_6^2\eta_7^2\eta_8}{\eta_1}\right|&\leq B. \label{H215}
\end{align}

\begin{lemma}(\cite[Lemma 8.4]{DE2})\label{L22}
We have
$$N_{U_2,H}(B) = \sum_{\substack{\boldsymbol{\eta}' \in \mathbb{Z}_1\times\cdots\times\mathbb{Z}_7}}
    \vartheta_1(\boldsymbol{\eta}')V_1(\boldsymbol{\eta}';B) + O\left(B(\log B)^2\right),$$
where
\begin{align}\label{E21}
V_1(\boldsymbol{\eta}';B)=\frac{1}{\eta_1}\int_{(\boldsymbol{\eta}', \eta_8)\in\mathcal{R}(B)}\mathrm{d}\eta_8
\end{align}
and
\begin{align}\label{E22}
\vartheta_1(\boldsymbol{\eta}')=\prod\limits_{p}\vartheta_{1,p}\left(I_p(\boldsymbol{\eta}')\right)
\end{align}
with $I_p(\boldsymbol{\eta}')=\left\{i\in\{1,\ldots,7\}:p\mid \eta_i\right\}$ and
\begin{align*}
\vartheta_{1,p}(I)=
\begin{cases}
     1,  &I=\emptyset, \{1\}, \{2\}, \{7\},\\
    1-\frac{1}{p},&I=\{4\}, \{5\}, \{6\}, \{1,3\}, \{2,3\}, \{3,4\}, \{4,6\},\{5,6\},\{5,7\},\\
    1-\frac{2}{p}, &I =\{3\},\\
    0, & \text{all other} \;I \subset \{1, \dots, 7\}.
\end{cases}
\end{align*}
\end{lemma}

\subsection{The second summation over $\eta_7$}
\begin{lemma}\label{L23}
We have
\begin{align*}
N_{U_2,H}(B) = \sum_{\substack{\boldsymbol{\eta} \in \mathbb{Z}_1\times\cdots\times\mathbb{Z}_6}}
   \mathcal{A}(\vartheta_1(\boldsymbol{\eta}'),\eta_7)V_{11}(\boldsymbol{\eta};B) + O\left(B(\log B)^4\log\log B\right),
\end{align*}
where $\mathcal{A}(\vartheta_1(\boldsymbol{\eta}'),\eta_7)$ is defined in Section 2 and
\begin{align}\label{V21}
V_{11}(\boldsymbol{\eta};B) =\frac{1}{\eta_1}\int_{(\boldsymbol{\eta}', \eta_8)\in\mathcal{R}(B)}\mathrm{d}\eta_7\mathrm{d}\eta_8.
\end{align}
\end{lemma}
\begin{proof}
Firstly, we rewrite the result of Lemma \ref{L22} as follow.
\begin{align*}
N_{U_2,H}(B) = \sum_{\substack{\boldsymbol{\eta} \in \mathbb{Z}_1\times\cdots\times\mathbb{Z}_6}}\sum_{\eta_7\in\mathbb{Z}_7}
    \vartheta(\eta_7)g(\eta_7) + O\left(B(\log B)^2\right),
\end{align*}
where\begin{align*}
\vartheta(\eta_7):= \vartheta_1(\boldsymbol{\eta},\eta_7) \qquad\text{and}\qquad g(\eta_7):= V_1(\boldsymbol{\eta},\eta_7;B).
\end{align*}

Since (\ref{E22}) and Definition \ref{D1}, we have $c=1$ and $A_p(0)=1$. Further, by Definition \ref{D2}, we yield that
\begin{align*}
\vartheta(\eta_7)\in \Theta_2(b,1,1,1)\qquad\text{and}\qquad b=\prod\limits_{p\mid \eta_1\eta_2\eta_3\eta_4\eta_6}p.
\end{align*}
Applying \cite[Lemma 5.1(4)]{DE2} to (\ref{H215}), we imply that
\begin{equation}\label{E24}
\begin{aligned}
g(\eta_7)&\ll\left(\frac{B}{\eta_1\eta_2|\eta_7|}\right)^\frac{1}{2}\\
&=\frac{B}{\boldsymbol{\eta}^{(1,1,1,1,1,1)}|\eta_7|}\left(\frac{B}{\boldsymbol{\eta}^{(2,2,3,2,0,1)}}\right)^{-1/4}
\left(\frac{B}{\boldsymbol{\eta}^{(0,0,1,2,4,3)}\eta_7^2}\right)^{-1/4}.
\end{aligned}
\end{equation}
In particular, by (\ref{H212}) and (\ref{H214}), we deduce that
\begin{align}\label{E25}
g(\eta_6)\ll\frac{B}{\boldsymbol{\eta}^{(1,1,1,1,1,1)}|\eta_7|}.
\end{align}

Let $t_1 := (\log B)^8$. Similar to Lemma \ref{L3}, we split the sum over $\eta_7$ into the two cases $|\eta_7| \leq t_1$ and $|\eta_7| > t_1$.
For the second case, by Lemma \ref{L11}, Remark \ref{R1} and (\ref{E25}), we get that
\begin{align*}
\sum_{|\eta_7| > t_1}\vartheta(\eta_7)g(\eta_7)=\mathcal{A}(\vartheta(\eta_7),\eta_7)\int\limits_{\substack{\eta_7\in\mathbb{Z}_7\\|\eta_7| > t_1}}g(\eta_7)\mathrm{d}\eta_7
+O\left(\frac{B\log B2^{\omega(\eta_1\eta_2\eta_3\eta_4\eta_6)}}{\boldsymbol{\eta}^{(1,1,1,1,1,1)}|\eta_7|}\right).
\end{align*}
When summing the error term over the remaining variables with $|\eta_7| > t_1$, we have
\begin{align*}
\ll \frac{B\log B}{|\eta_7|}\sum_{\substack{\boldsymbol{\eta}}}\frac{2^{\omega(\eta_1\eta_2\eta_3\eta_4\eta_6)}}{\boldsymbol{\eta}^{(1,1,1,1,1,1)}}
\ll t_1^{-1}B(\log B)^{12}\ll B(\log B)^4.
\end{align*}

Now, we start with the first case. Since $|\eta_7|\leq t_1$ and $0\leq\vartheta(\eta_7)\leq 1$, by (\ref{E24}), we yield an upper bound
\begin{equation}\label{E26}
\begin{aligned}
\sum_{\substack{\boldsymbol{\eta} \in \mathbb{Z}_1\times\cdots\times\mathbb{Z}_6}}&\sum_{\substack{\eta_7\in\mathbb{Z}_7\\|\eta_7| \leq t_1}}\vartheta(\eta_7)g(\eta_7)\\
&\ll \sum_{\substack{\boldsymbol{\eta}'\\ |\eta_7| \leq t_1}}\frac{B}{\boldsymbol{\eta}^{(1,1,1,1,1,1)}|\eta_7|}\left(\frac{B}{\boldsymbol{\eta}^{(2,2,3,2,0,1)}}\right)^{-1/4}
\left(\frac{B}{\boldsymbol{\eta}^{(0,0,1,2,4,3)}\eta_7^2}\right)^{-1/4}\\
&\ll \sum_{\substack{\eta_2,\eta_3,\eta_4,\eta_5,\eta_6,\eta_7\\ |\eta_7| \leq t_1}}\frac{B}{\eta_2\eta_3\eta_4\eta_5\eta_6|\eta_7|}\left(\frac{B}{\boldsymbol{\eta}^{(0,0,1,2,4,3)}\eta_7^2}\right)^{-1/4}\\
&\ll  \sum_{\substack{\eta_2,\eta_4,\eta_5,\eta_6,\eta_7\\ |\eta_7| \leq t_1}}\frac{B}{\eta_2\eta_4\eta_5\eta_6|\eta_7|}\ll B(\log B)^4\log t_1\\
&\ll B(\log B)^4\log\log B,
\end{aligned}
\end{equation}
where we have used the conditions (\ref{H212}) and (\ref{H214}).

In order to prove this lemma, it suffices to obtain the following upper bound.
\begin{align*}
&\sum_{\substack{\boldsymbol{\eta} \in \mathbb{Z}_1\times\cdots\times\mathbb{Z}_6}}\mathcal{A}(\vartheta(\eta_7),\eta_7)\int\limits_{\substack{\eta_7\in\mathbb{Z}_7\\|\eta_7| \leq t_1}}g(\eta_7)\mathrm{d}\eta_7\\
&\ll\sum_{\substack{\boldsymbol{\eta}}}\int\limits_{|\eta_7| \leq t_1}\frac{B}{\boldsymbol{\eta}^{(1,1,1,1,1,1)}|\eta_7|}\left(\frac{B}{\boldsymbol{\eta}^{(2,2,3,2,0,1)}}\right)^{-1/4}
\left(\frac{B}{\boldsymbol{\eta}^{(0,0,1,2,4,3)}\eta_7^2}\right)^{-1/4}\mathrm{d}\eta_7\\
&\ll\sum_{\substack{\eta_2,\eta_4,\eta_5,\eta_6}}\frac{B}{\eta_2\eta_4\eta_5\eta_6}\int\limits_{|\eta_7| \leq t_1}|\eta_7|^{-1}\mathrm{d}\eta_7
\ll B(\log B)^4\log\log B,
\end{align*}
where we have used (\ref{E26}), $0\leq \mathcal{A}(\vartheta(\eta_7),\eta_7)\leq 1$.

We complete this lemma, since
\begin{align*}
\mathcal{A}(\vartheta(\eta_7),\eta_7)=\mathcal{A}(\vartheta_1(\boldsymbol{\eta}'),\eta_7)\qquad\text{and}\qquad\int\limits_{\eta_7\in\mathbb{Z}_7} g(\eta_7)\mathrm{d}\eta_7=V_{11}(\boldsymbol{\eta};B).
\end{align*}
\end{proof}

\subsection{The remaining summations}
\begin{lemma}\label{L24}
We have
\begin{align*}
N_{U_2,H}(B)=\left(\prod\limits_{p}\omega_p\right)V_0(B)+O\left(B(\log B)^4\log\log B\right),
\end{align*}
where
\begin{align*}
\omega_p:=\left(1-\frac{1}{p}\right)^6\left(1+\frac{6}{p}+\frac{1}{p^2}\right)
\end{align*} and
\begin{align*}
V_0(B): =\int_{(\boldsymbol{\eta}', \eta_8)\in\mathcal{R}(B)}\frac{1}{\eta_1}\mathrm{d}\eta_1\cdots\mathrm{d}\eta_8.
\end{align*}
\end{lemma}
\begin{proof}
Firstly, we rewrite the result of Lemma \ref{L23} as follow.
\begin{align}\label{V22}
N_{U_2,H}(B) =\sum_{\substack{\boldsymbol{\eta} \in \mathbb{Z}_1\times\cdots\times\mathbb{Z}_6}}
   \vartheta_6(\boldsymbol{\eta})V_6(\boldsymbol{\eta};B)+O\left(B(\log B)^4\log\log B\right),
\end{align}
where\begin{align*}
\vartheta_6(\boldsymbol{\eta}):= \mathcal{A}(\vartheta_1(\boldsymbol{\eta}'),\eta_7)\qquad\text{and}\qquad
V_6(\boldsymbol{\eta};B):= V_{11}(\boldsymbol{\eta};B).
\end{align*}

Since (\ref{E22}), we have $\vartheta_1(\boldsymbol{\eta}')\in\Theta'_{4,7}(2)$ \cite[Definition 7.8]{DE2}. By \cite[Fig.2]{DE2}, we get $\vartheta_1(\boldsymbol{\eta}')\in\Theta_{2,7}(C)$ \cite[Definition 4.2]{DE2} for some $C\in\mathbb{Z}_{\geq0}$. Furthermore, we have
$\vartheta_6(\boldsymbol{\eta})=\mathcal{A}(\vartheta_1(\boldsymbol{\eta}'),\eta_7)\in\Theta_{2,6}(C)$. Applying \cite[Lemma 5.1(5)]{DE2} to (\ref{H215}),
we have
\begin{align*}
V_6(\boldsymbol{\eta};B)\ll\frac{B^{2/3}}{(\eta_1\eta_2\eta_4\eta_5^3\eta_6^2)^{1/3}}
=\frac{B}{\boldsymbol{\eta}^{(1,1,1,1,1,1)}}\left(\frac{B}{\boldsymbol{\eta}^{(2,2,3,2,0,1)}}\right)^{-1/3}.
\end{align*}
Now, the conditions of \cite[Proposition 4.3]{DE2} are satisfied, and we apply it (see (\ref{A})) to the main term of (\ref{V22}) with $(r,s)=(5,1)$, we imply that
\begin{align*}
\sum_{\substack{\boldsymbol{\eta} \in \mathbb{Z}_1\times\cdots\times\mathbb{Z}_6}}
   \vartheta_6(\boldsymbol{\eta})V_6(\boldsymbol{\eta};B)=\mathcal{A}( \vartheta_6(\boldsymbol{\eta}),\boldsymbol{\eta}) \int\limits_{\boldsymbol{\eta}}V_6(\boldsymbol{\eta};B)\mathrm{d}\boldsymbol{\eta}+O\left(B(\log B)^4\log\log B\right).
\end{align*}
By (\ref{V21}), we have
\begin{align*}
\int\limits_{\boldsymbol{\eta}}V_6(\boldsymbol{\eta};B)\mathrm{d}\boldsymbol{\eta}=V_0(B).
\end{align*}

Since $\mathcal{A}( \vartheta_6(\boldsymbol{\eta}),\boldsymbol{\eta})=\mathcal{A}(\vartheta_1(\boldsymbol{\eta}'),\boldsymbol{\eta}')$ and
$\vartheta_1(\boldsymbol{\eta}')\in\Theta'_{4,7}(2)$, we get that the conditions of \cite[Corollary 7.10]{DE2} are satisfied.
From \cite[Corollary 7.10]{DE2} and (\ref{E22}), we have
\begin{align}\label{W}
\mathcal{A}( \vartheta_6(\boldsymbol{\eta}),\boldsymbol{\eta})=\prod\limits_{p}\omega_p,
\end{align}
where
\begin{align*}
\omega_p=&\left(1-\frac{1}{p}\right)^7+3\left(1-\frac{1}{p}\right)^6\left(\frac{1}{p}\right)
+3\left(1-\frac{1}{p}\right)^6\left(\frac{1}{p}\right)\left(1-\frac{1}{p}\right)\\
&+6\left(1-\frac{1}{p}\right)^5\left(\frac{1}{p}\right)^2\left(1-\frac{1}{p}\right)
+\left(1-\frac{1}{p}\right)^6\left(\frac{1}{p}\right)\left(1-\frac{2}{p}\right)\\
=&\left(1-\frac{1}{p}\right)^6\left(1+\frac{6}{p}+\frac{1}{p^2}\right).
\end{align*}
\end{proof}

\subsection{The final proof}
$\newline$ Now, we will show that the volume of this region grows asymptotically as predicted by Manin and Peyre, i.e. step (3) in Section 1. Since the proof is the same to \cite[Section 8.3]{DE2}, we only give the outline.

Define
\begin{equation*}
  \begin{split}
    &\mathcal{R}'_1(B)=\{(\eta_1, \dots, \eta_5) \in J'_1\times \dots \times J'_5 \mid
    \eta_1^2\eta_2^2\eta_3^3\eta_4^2 \le B,\ \eta_1^3\eta_2^3\eta_3^4\eta_4^2\eta_5^{-2}
    \ge B\},\\
    &\mathcal{R}'_2(B)=(\eta_1, \dots, \eta_5;B) = \{(\eta_6,\eta_7,\eta_8) \in J'_6\times J'_7
    \times J'_8 \mid (\ref{H211})-(\ref{H215})\},\\
    &\mathcal{R}'(B)= \{(\eta_1, \dots,\eta_8) \in \mathbb{R}^8 \mid (\eta_1, \dots, \eta_5) \in
   \mathcal{R}'_1(B), (\eta_6,\eta_7,\eta_8) \in \mathcal{R}'_2(B)\}
  \end{split}
\end{equation*}
and
\begin{equation*}
  V_0'(B) = \int_{(\boldsymbol{\eta}', \eta_8)\in \mathcal{R}'(B)} \eta_1^{-1} \mathrm{d}\eta_1\cdots\mathrm{d}\eta_8.
\end{equation*}

From \cite[Lemma 8.6]{DE2} and \cite[Lemma 8.7]{DE2}, we have
\begin{equation}\label{E27}
 V_0(B) = V_0'(B) + O(B(\log B)^4),
\end{equation}
where
\begin{equation}\label{E28}
  V_0'(B) =\alpha(\widetilde{S}_2)\prod\limits_{p}\left(1-\frac{1}{p}\right)^6\left(1+\frac{6}{p}+\frac{1}{p^2}\right)\omega_\infty(\widetilde{S}_2)B(\log B)^5
\end{equation}
with $\alpha(\widetilde{S}_2)=\frac{1}{8640}$ from \cite{DE2}.

Theorem 1.2 follows from (\ref{E27})-(\ref{E28}) and Lemma \ref{L24}.

\section{Proof for the quartic del Pezzo surface $S_3$ of type $\mathbf{A}_4$}
Let $S_3$, $\widetilde{S}_3$, $U_3$ and $N_{U_3,H}(B)$ be defined in Section 1. The proof of Theorem 1.3 is similar to Theorem 1.1. However, it is slightly difficult than $S_1$, since we need to take extra efforts to deal with a quadratic congruence condition (see (\ref{Quadratic})). Our proof follows the idea of \cite[Section 3]{DEFR2}.

\subsection{Passage to a universal torsor}
$\newline$We will introduce new variables $\eta_1,\ldots,\eta_9$ and use the notation
\begin{align*}
\boldsymbol{\eta}:=(\eta_1,\ldots,\eta_6),\quad\boldsymbol{\eta}':=(\eta_1,\ldots,\eta_7),\quad\boldsymbol{\eta}'':=(\eta_1,\ldots,\eta_9)
\quad\boldsymbol{\eta}''':=(\eta_1,\ldots,\eta_5,\eta_7)
\end{align*}
and, for $(k_1,\ldots,k_6)\in\mathbb{Q}^6$,
\begin{align*}
\boldsymbol{\eta}^{(k_1,k_2,k_3,k_4,k_5,k_6)}:=\eta_1^{k_1}\eta_2^{k_2}\eta_3^{k_3}\eta_4^{k_4}\eta_5^{k_5}\eta_6^{k_6}.
\end{align*}
For $i=1,\ldots,9$, let
\begin{equation*}
  (\mathbb{Z}_i, J_i, J_i')=
  \begin{cases}
    (\mathbb{Z}_{>0}, \mathbb{R}_{\ge 1}, \mathbb{R}_{\ge 1}), & i\in \{1, 2\},\\
   (\mathbb{Z}_{>0}, \mathbb{R}_{\ge 1}, \mathbb{R}), & i=3,\\
    (\mathbb{Z}_{>0}, \mathbb{R}_{\ge 1}, \mathbb{R}_{\ge 1}), & i\in \{4,5,6\},\\
     (\mathbb{Z}_{\neq 0}, \mathbb{R}_{\le -1} \cup \mathbb{R}_{\ge 1}, \mathbb{R}), & i=7,\\
    (\mathbb{Z}, \mathbb{R},\mathbb{R}), & i \in \{8, 9\}.
  \end{cases}
\end{equation*}
In this section, we estimate summations over $\eta_i\in\mathbb{Z}_i$ by integrations over $\eta_i\in J_i$, which we enlarge to $\eta_i\in J_i'$ in Lemma \ref{L38}. Next, our purpose is to establish a bijection between the rational points on the open subset $U_3$, and the integral points on the universal torsor above $\widetilde{S}_3$, which are subject to a number of coprimality conditions.

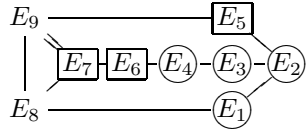
\begin{figure}[ht]
  \centering
  \[\xymatrix@R=0.05in @C=0.05in{E_9 \ar@{-}[rrrr] \ar@1{=}[dr] \ar@1{-}[dd] &&&&  \li{5} \ar@{-}[dr]\\
      & \li{7} \ar@{-}[r] & \li{6} \ar@{-}[r] & \ex{4}\ar@{-}[r] & \ex{3}\ar@{-}[r] & \ex{2}\\
      E_8 \ar@{-}[rrrr] \ar@{-}[ur] &&&& \ex{1} \ar@{-}[ur]}\]
  \caption{Configuration of curves on $\widetilde{S}_3$.}
  \label{F3}
\end{figure}

Firstly, we must begin by collecting together some useful information about the geometric structure of $S_3$ from \cite[Section 3.4]{DE1}. The singularity $\mathbf{A}_4$ in $(0 : 0 : 0 : 0 : 1)$ gives exceptional divisors $E_1, E_2, E_3, E_4$. Let $E_5, E_6, E_7$ resp. $E_8, E_9$ on $\widetilde{S}_3$
be the strict transforms under $\pi:\widetilde{S}_3\to S_3$ of the three lines $E_5''=\{x_0=x_2=x_3=0\}$, $E_6''=\{x_0=x_1=x_3=0\}$, $E_7''=\{x_1=x_3=x_4=0\}$ resp. the curves $E_8''=\{x_1=x_2=x_0x_4+x_3^2=0\}$, $E_9''=\{x_4=x_0x_1-x_2x_3=x_1x_2+x_3^2=x_0x_3+x_2^2=0\}$ on $S_3$. The configuration of these divisors $E_1,\ldots,E_9$ on $\widetilde{S}_3$ is described by Figure \ref{F3} with circles and squares marking classes of $(-2)$- and $(-1)$-curves.

\begin{lemma}\label{L31}
 We have
$$N_{U_3,H}(B) = \#\mathcal{T}_3(B),$$
where $\mathcal{T}_3(B)$ is the set of $\boldsymbol{\eta}'' \in \mathbb{Z}_1\times\cdots\times\mathbb{Z}_9$
such that
\begin{equation}\label{U31}
\eta_5\eta_9+\eta_1\eta_8^2+\eta_3\eta_4^2\eta_6^3\eta_7=0,
\end{equation}
holds, with
\begin{equation}\label{H31}
\max\left\{\begin{aligned}
&|\boldsymbol{\eta}^{(2,4,3,2,3,1)}|,|\boldsymbol{\eta}^{(1,1,1,1,0,1)}\eta_7\eta_8|, \\ &|\boldsymbol{\eta}^{(2,3,2,1,2,0)}\eta_8|,|\boldsymbol{\eta}^{(1,2,2,2,1,2)}\eta_7|,|\eta_7\eta_9|
\end{aligned}\right\} \le B,
\end{equation}
\begin{align}\label{C31}
\boldsymbol{\eta}''\; \text{fulfills coprimality conditions as in Figure \ref{F3}}.
\end{align}
\end{lemma}
\begin{proof}
 It is based on the construction of the minimal desingularization $\pi: \widetilde{S}_3 \to S_3$ by the following blow-ups
 $\rho: \widetilde{S}_3 \to \mathbb{P}^2$ in five points. Starting with the curves
\begin{equation*}
 E^{(0)}_3:= \{y_0=0\},\;E^{(0)}_8:=\{y_1=0\},\; E^{(0)}_7:=\{y_2 = 0\},\; E^{(0)}_9:=\{-y_0y_2-y_1^2= 0\}
\end{equation*}
in $\mathbb{P}^2$:
  \begin{enumerate}
 \item blow up $E_3^{(0)}\cap E_8^{(0)} \cap E_9^{(0)}$, giving $E_1^{(1)}$,
  \item blow up $E_1^{(1)}\cap E_3^{(1)} \cap E_9^{(1)}$, giving $E_2^{(2)}$,
  \item blow up $E_2^{(2)} \cap E_9^{(2)}$, giving $E_5^{(3)}$,
  \item blow up $E_3^{(3)} \cap E_7^{(3)}$, giving $E_4^{(4)}$,
  \item blow up $E_4^{(4)} \cap E_7^{(4)}$, giving $E_6^{(5)}$.
  \end{enumerate}

We may choose $\rho(E_3), \rho(E_8), \rho(E_7)$ as the coordinate lines in $\mathbb{P}^2 = \{(\eta_3: \eta_8: \eta_7)\}$. Then $\rho(E_9)$
is the quadric $\eta_3\eta_7+\eta_8^2=0$. The morphisms $\rho, \pi$ and the projection
\begin{equation*}
  \begin{array}[h]{cccc}
    \phi: & S_3 & \dashrightarrow & \mathbb{P}^2,\\
    & \mathbf{x} & \mapsto & (x_0:x_2:x_3)
  \end{array}
\end{equation*}
form a following commutative diagram of rational maps between $S_3, \widetilde{S}_3$ and $\mathbb{P}^2$.
\begin{equation*}
  \xymatrix{\widetilde{S}_3 \ar@{->}^\rho[dr]\ar@{->}^\pi[d] & \\
    S_3\ar@{-->}^\phi[r] & \mathbb{P}^2}
\end{equation*}
The inverse map of $\phi$ is
\begin{equation}\label{P31}
  \begin{array}[h]{cccc}
    \phi': & \mathbb{P}^2&\dashrightarrow &S_3,\\
    &(\eta_3: \eta_8: \eta_7)&\mapsto&(\eta_3^3:\eta_3\eta_7\eta_8:\eta_3^2\eta_8:\eta_3^2\eta_7:\eta_7\eta_9)
  \end{array}
\end{equation}
where $\eta_9=-\eta_3\eta_7-\eta_8^2$. The maps $\phi, \phi^\prime$ give a bijection between the complement $U_3$ of the lines on $S_3$ and $\{(\eta_3: \eta_8: \eta_7) \in \mathbb{P}^2\mid \gcd(\eta_3,\eta_7,\eta_8)=1\}$, and furthermore, induces a bijection between $U_3(\mathbb{Q})$ and the integral points
\begin{equation*}
  \{(\eta_3, \eta_7, \eta_8, \eta_9) \in \mathbb{Z}_{>0}\times\mathbb{Z}_{\neq 0}\times \mathbb{Z}^2 \mid
  \gcd(\eta_3, \eta_8, \eta_9)=1,\ \eta_3\eta_7+\eta_8^2+\eta_9=0\}.
\end{equation*}
Since the above five blow-ups (1)-(5), we introduce the following further variables:
\begin{equation*}
  \begin{array}[h]{lll}
    \eta_1=\gcd(\eta_3,\eta_8,\eta_9),&\eta_2=\gcd(\eta_1,\eta_3,\eta_9),&
   \eta_5=\gcd(\eta_2,\eta_9),\\
    \eta_4=\gcd(\eta_3,\eta_7),&\eta_6=\gcd(\eta_4,\eta_7).
  \end{array}
\end{equation*}

After the similar procedure to 1)-4) in Lemma \ref{L1}, we get (\ref{U31}), (\ref{C31}) and
\begin{align*}
\phi'': &(\eta_1 : \eta_2 : \eta_3 : \eta_4 : \eta_5 : \eta_6 : \eta_7 : \eta_8 : \eta_9)\mapsto\\
&(\boldsymbol{\eta}^{(2,4,3,2,3,1)}:\boldsymbol{\eta}^{(1,1,1,1,0,1)}\eta_7\eta_8:\boldsymbol{\eta}^{(2,3,2,1,2,0)}\eta_8:
\boldsymbol{\eta}^{(1,2,2,2,1,2)}\eta_7:\eta_7\eta_9).
\end{align*}
Therefore, $H(\phi''(\boldsymbol{\eta}'')) \le B$ is equivalent to (\ref{H31}).
\end{proof}

\subsection{The first summation over $\eta_8$}
$\newline$Similar to Table \ref{table1}, the Table \ref{table2} also provides a dictionary between the notation of \cite[Section 2]{DE2} and the present situation on $S_3$.
\begin{table}[!ht]
\centering
\caption{Dictionary for applying \cite[Proposition 2.4]{DE2}}
\begin{tabular}{|c|c|c|c|} \hline
$(r,s,t)$ & $(3,1,1)$& $ \delta$ & $\eta_2$\\ \hline
($\alpha_0$;$\alpha_1,\ldots,\alpha_r$) & $(\eta_7;\eta_3,\eta_4,\eta_6)$ &
$(a_0;a_1,\ldots,a_r)$ & $(1;1,2,3)$ \\ \hline
($\beta_0$;$\beta_1,\ldots,\beta_s$) & $(\eta_8;\eta_1)$ &
$(b_0;b_1,\ldots,b_s)$ & $(2;1)$ \\ \hline
($\gamma_0$;$\gamma_1,\ldots,\gamma_t$) & $(\eta_9;\eta_5)$ &
$(c_1,\ldots,c_t)$ & $(1)$ \\ \hline
$\Pi(\boldsymbol{\alpha})$ & $\eta_3\eta_4^2\eta_6^3$ & $\Pi'(\delta,\boldsymbol{\alpha})$ & $\eta_2\eta_3\eta_4$\\ \hline
$\Pi(\boldsymbol{\beta})$ & $\eta_1$ & $\Pi'(\delta,\boldsymbol{\beta})$ & $\eta_2$\\ \hline
$\Pi(\boldsymbol{\gamma})$ & $\eta_5$ & $\Pi'(\delta,\boldsymbol{\gamma})$ & $\eta_2$\\
\hline
\end{tabular}
\label{table2}
\end{table}

Using (\ref{U31}) to eliminate $\eta_9$, from the height condition (\ref{H31}), we define $\mathcal{R}(B)$ to be the set of all $(\boldsymbol{\eta}',\eta_8) \in J_1\times\cdots \times J_8$ and
\begin{align}
\left|\boldsymbol{\eta}^{(2,4,3,2,3,1)}\right|&\leq B, \label{H311}\\
\left|\boldsymbol{\eta}^{(1,1,1,1,0,1)}\eta_7\eta_8\right|&\leq B, \label{H312}\\
\left|\boldsymbol{\eta}^{(2,3,2,1,2,0)}\eta_8\right|&\leq B, \label{H313}\\
\left|\boldsymbol{\eta}^{(1,2,2,2,1,2)}\eta_7\right|&\leq B, \label{H314}\\
\left|\frac{\eta_1\eta_7\eta_8^2+\eta_3\eta_4^2\eta_6^3\eta_7^2}{\eta_5}\right|&\leq B. \label{H315}
\end{align}

\begin{lemma}\label{L32}
We have
$$N_{U_3,H}(B) = \sum_{\substack{\boldsymbol{\eta}' \in \mathbb{Z}_1\times\cdots\times\mathbb{Z}_7}}
    \vartheta_1(\boldsymbol{\eta}')V_1(\boldsymbol{\eta}';B) + O\left(B(\log B)^3\right),$$
where
\begin{align}\label{E31}
V_1(\boldsymbol{\eta}';B)=\frac{1}{\eta_5}\int_{(\boldsymbol{\eta}', \eta_8)\in\mathcal{R}(B)}\mathrm{d}\eta_8
\end{align}
and
\begin{equation}\label{E32}
\vartheta_1(\boldsymbol{\eta}')=\sum_{\substack{k|\eta_2\\ \gcd(k,\eta_1\eta_3)}}
\frac{\mu(k)}{k}\widetilde{\vartheta}_1(\boldsymbol{\eta}',k)
\sum_{\substack{\varrho^2\equiv-\eta_3\eta_6\eta_7/\eta_1\pmod{k\eta_5}\\ 1\leq\varrho\leq k\eta_5\\ \gcd(\varrho,k\eta_5)}}1
\end{equation}
with
\begin{equation}\label{E32!}
\widetilde{\vartheta}_1(\boldsymbol{\eta}',k)=\frac{\phi^\ast(\eta_2\eta_3\eta_4\eta_6)}{\phi^\ast(\gcd(\eta_2,k\eta_5))}.
\end{equation}
\end{lemma}
\begin{proof}
From Table \ref{table2} and \cite[Proposition 2.4]{DE2}, we have
\begin{align*}
N_{U_3,H}(B) =\sum_{\substack{\boldsymbol{\eta}' \in\mathbb{Z}_1\times\cdots\times\mathbb{Z}_7}}
     \Big(\vartheta_1(\boldsymbol{\eta}')V_1(\boldsymbol{\eta}';B) +R_1(\boldsymbol{\eta}';B)\Big),
\end{align*}
where $V_1(\boldsymbol{\eta}';B)$ is defined as (\ref{E31}),
\begin{align}\label{Quadratic}
\vartheta_1(\boldsymbol{\eta}')=\sum\limits_{\substack{k\mid \eta_2\\ \gcd(k,\eta_1\eta_3\eta_4\eta_6\eta_7)=1}}
 \frac{\mu(k)\phi^\ast(\eta_2\eta_3\eta_4\eta_6)}{k\phi^\ast(\gcd(\eta_2,k\eta_5))}
  \sum_{\substack{\varrho^2\equiv-\eta_3\eta_6\eta_7/\eta_1\pmod{k\eta_5}\\ 1\leq\varrho\leq k\eta_5\\ \gcd(\varrho,k\eta_5)}}1
\end{align}
and\begin{align*}
R_1(\boldsymbol{\eta}';B)\ll 2^{\omega(\eta_2)+\omega(\eta_2\eta_3\eta_4\eta_6)+\omega(\eta_2\eta_5)}.
\end{align*}

In order to decide the main term, it suffices to compute $\vartheta_1(\boldsymbol{\eta}')$. Since $k|\eta_2$ and (\ref{C31}), we have $\gcd(k,\eta_4\eta_6\eta_7)=1$. Thus, (\ref{Quadratic}) is equivalent to (\ref{E32}).

Next, we will utilize (\ref{H314}) to deal with the error term.
\begin{align*}
\sum_{\substack{\boldsymbol{\eta}'}}R_1(\boldsymbol{\eta}';B)
&\ll \sum_{\substack{\boldsymbol{\eta}'}}2^{\omega(\eta_2)+\omega(\eta_2\eta_3\eta_4\eta_6)+\omega(\eta_2\eta_5)}\\
&\ll B\sum_{\substack{\boldsymbol{\eta}}}\frac{2^{\omega(\eta_2)+\omega(\eta_2\eta_3\eta_4\eta_6)+\omega(\eta_2\eta_5)}}{\boldsymbol{\eta}^{(1,2,2,2,1,2)}}\\
&\ll B(\log B)^3.
\end{align*}
\end{proof}

If applying Lemma \ref{L11} to the main term in Lemma \ref{L32} directly, we can not yield the desired error term. Thus, we will divide it into the two cases $\eta_6>|\eta_7|$ and $\eta_6\leq|\eta_7|$. Let $N_6(B)$ be the main term in Lemma \ref{L32} with the additional condition $\eta_6>|\eta_7|$ on the $\boldsymbol{\eta}'$. Similarly, $N_7(B)$ be the main term in Lemma \ref{L32} with the additional condition $\eta_6\leq|\eta_7|$ on the $\boldsymbol{\eta}'$.
Hence, we have
\begin{align}\label{N}
N_{U_3,H}(B)=N_6(B)+N_7(B)+ O\left(B(\log B)^3\right).
\end{align}

Moreover, let $\vartheta_0(\boldsymbol{\eta}'):=\prod\limits_{p}\vartheta_{0,p}\left(I_p(\boldsymbol{\eta}')\right)$
with $I_p(\boldsymbol{\eta}')=\left\{i\in\{1,\ldots,7\}:p\mid \eta_i\right\}$ and
\begin{align*}
\vartheta_{0,p}(I)=
\begin{cases}
     1,  &I=\emptyset, \{1\}, \{2\}, \{3\}, \{4\}, \{5\}, \{6\}, \{7\},\\
         &\text{or}\;I=\{1,2\}, \{2,3\}, \{2,5\}, \{3,4\}, \{4,6\},\{6,7\}\\
    0, & \text{all other} \;I \subset \{1, \dots, 7\}.
\end{cases}
\end{align*}
Then
\begin{align}\label{Special2}
\vartheta_0(\boldsymbol{\eta}')=1
\end{align}
if and only if $\boldsymbol{\eta}'$ fulfills coprimality conditions from (\ref{C31}).

Define\begin{align}\label{Quadratic1}
\vartheta'_1(\boldsymbol{\eta}'):=\sum\limits_{\substack{k\mid \eta_2\\ \gcd(k,\eta_1\eta_3)=1}}
 \frac{\mu(k)\phi^\ast(\eta_2\eta_3\eta_4\eta_6)}{k\phi^\ast(\gcd(\eta_2,k\eta_5))}.
\end{align}

\subsection{The asymptotic formula on $N_6(B)$}
\subsubsection{The second summation over $\eta_6$}
\begin{lemma}\label{L33}
We have
\begin{align*}
N_6(B)= \sum_{\substack{\boldsymbol{\eta}''' \in \mathbb{Z}_1\times\cdots\times\mathbb{Z}_5\times\mathbb{Z}_7}}
   \mathcal{A}(\vartheta'_1(\boldsymbol{\eta}'),\eta_6)V_{16}(\boldsymbol{\eta}''';B)+O\left(B(\log B)^4\right),
\end{align*}
where $\mathcal{A}(\vartheta'_1(\boldsymbol{\eta}'),\eta_6)$ is defined in Section 2 and
\begin{align}\label{V31}
V_{16}(\boldsymbol{\eta}''';B) =\frac{1}{\eta_5}\int\limits_{\substack{(\boldsymbol{\eta}', \eta_8)\in\mathcal{R}(B)\\ \eta_6>|\eta_7|}}
\mathrm{d}\eta_6\mathrm{d}\eta_8.
\end{align}
\end{lemma}
\begin{proof}
Firstly, by Lemma \ref{L32}, we have
\begin{align}\label{E33}
N_6(B) = \sum_{\substack{\boldsymbol{\eta}''' \in \mathbb{Z}_1\times\cdots\times\mathbb{Z}_5\times\mathbb{Z}_7}}\sum_{\substack{k|\eta_2\\ \gcd(k,\eta_1\eta_3)}}
\frac{\mu(k)}{k}\Sigma,
\end{align}
where
\begin{align*}
\Sigma:=\sum_{\substack{\eta_6 \in \mathbb{Z}_6\\\eta_6>|\eta_7|}}\vartheta(\eta_6)\sum_{\substack{\varrho^2\equiv-\eta_3\eta_6\eta_7/\eta_1\pmod{k\eta_5}\\ 1\leq\varrho\leq k\eta_5\\ \gcd(\varrho,k\eta_5)}}g(\eta_6)
\end{align*}
with
\begin{align*}
\vartheta(\eta_6):= \widetilde{\vartheta}_1(\boldsymbol{\eta}',k) \qquad\text{and}\qquad g(\eta_6):= V_1(\boldsymbol{\eta}''',\eta_6;B).
\end{align*}

From (\ref{E32!}), (\ref{Special2}) and Definitions \ref{D1}-\ref{D2}, we yield immediately that
$c=1$, $A_p(0)=1$,
\begin{align*}
\vartheta(\eta_6)\in \Theta_2(b,1,1,1)\qquad\text{and}\qquad b=\prod\limits_{p\mid \eta_1\eta_2\eta_3\eta_5}p.
\end{align*}
Moreover, applying \cite[Lemma 5.1(1)]{DE2} to (\ref{H315}), we imply that
\begin{align}\label{E34}
g(\eta_6)\ll\left(\frac{B}{\eta_1\eta_5|\eta_7|}\right)^{1/2}.
\end{align}

From $\Sigma$, we learn that the sum over $\varrho$ is quadratic and is not 1. By Lemma \ref{L11}, we get that
\begin{align*}
\Sigma=\phi^\ast(k\eta_5)\mathcal{A}(\vartheta(\eta_6),\eta_6,k\eta_5)&\int\limits_{\substack{\eta_6\in\mathbb{Z}_6\\ \eta_6>|\eta_7|}}g(\eta_6)\mathrm{d}\eta_6\\
&+O\left(kB^{1/2}\log B\frac{\eta_5^{1/2}2^{\omega(\eta_1\eta_2\eta_3\eta_5)}}{\eta_1^{1/2}|\eta_7|^{1/2}}\right).
\end{align*}

Since (\ref{H311}) and $\eta_6>|\eta_7|$, we have
\begin{align}\label{E35}
\boldsymbol{\eta}^{(2,4,3,2,3,0)}|\eta_7|\leq B.
\end{align}
When summing the error term over $k$ and the remaining variables, we have
\begin{align*}
&\ll B^{1/2}\log B\sum_{\substack{\boldsymbol{\eta}'''\\(\ref{E35})}}
\frac{\eta_5^{1/2}2^{\omega(\eta_2)+\omega(\eta_1\eta_2\eta_3\eta_5)}}{\eta_1^{1/2}|\eta_7|^{1/2}}\\
&\ll B\log B\sum_{\substack{\eta_1,\cdots,\eta_5}}\frac{2^{\omega(\eta_2)+\omega(\eta_1\eta_2\eta_3\eta_5)}}{\eta_1^{3/2}\eta_2^2\eta_3^{3/2}\eta_4\eta_5}\\
&\ll B(\log B)^4.
\end{align*}
Now, we have
\begin{align*}
N_6(B) = \sum_{\substack{\boldsymbol{\eta}''' \in \mathbb{Z}_1\times\cdots\times\mathbb{Z}_5\times\mathbb{Z}_7}}\sum_{\substack{k|\eta_2\\ \gcd(k,\eta_1\eta_3)}}
\frac{\mu(k)}{k}\phi^\ast(k\eta_5)&\mathcal{A}(\vartheta(\eta_6),\eta_6,k\eta_5)V_{16}(\boldsymbol{\eta}''';B)\\
&+O\left(B(\log B)^4\right).
\end{align*}

In order to complete this lemma, it suffices to prove
\begin{align}\label{E36}
\sum_{\substack{k|\eta_2\\ \gcd(k,\eta_1\eta_3)}}\frac{\mu(k)}{k}\phi^\ast(k\eta_5)&\mathcal{A}(\vartheta(\eta_6),\eta_6,k\eta_5)
= \mathcal{A}(\vartheta'_1(\boldsymbol{\eta}'),\eta_6).
\end{align}
Next, we will use the idea from \cite[Lemma 6.3]{DEFR1} to prove (\ref{E36}). From (\ref{E32!}) and (\ref{Quadratic1}), it is enough to prove that
\begin{align}\label{E37}
\phi^\ast(k\eta_5)\mathcal{A}(\vartheta(\eta_6),\eta_6,k\eta_1)
= \mathcal{A}(\vartheta(\eta_6),\eta_6).
\end{align}
If $\vartheta(\eta_6)$ is the zero function, this is clearly true. If not, write $\vartheta(\eta_6)=\prod\limits_{p^\nu\|\eta_6}A_p(\nu)$ with $A_p(\nu)=A_p(\nu+1)$
for $\nu\geq 1$. By Definition \ref{D2}(3),
\begin{align*}
\phi^\ast(k\eta_5)\mathcal{A}(\vartheta(\eta_6),\eta_6,k\eta_1)
=\prod\limits_{p\nmid k\eta_5}\left((1-\frac{1}{p})A_p(0)+\frac{1}{p}A_p(1)\right)\prod\limits_{p| k\eta_5}\left(1-\frac{1}{p}\right)A_p(0)
\end{align*}
and
\begin{align*}
\mathcal{A}(\vartheta(\eta_6),\eta_6)=\prod\limits_{p}\left((1-\frac{1}{p})A_p(0)+\frac{1}{p}A_p(1)\right).
\end{align*}
Whenever $p|\gcd(\eta_6,k\eta_5)$, we have $\vartheta(\eta_6)=\widetilde{\vartheta}_1(\boldsymbol{\eta}',k)=0$. Since $\vartheta(\eta_6)$
is not identically zero, this implies $A_p(1)=0$ for all $p|k\eta_5$. This completes (\ref{E37}).
\end{proof}

\subsubsection{The remaining summations}
\begin{lemma}\label{L34}
We have
\begin{align*}
N_6(B)=\left(\prod\limits_{p}\omega_p\right)V_{60}(B)+O\left(B(\log B)^4\log\log B\right),
\end{align*}
where
\begin{align*}
\omega_p:=\left(1-\frac{1}{p}\right)^6\left(1+\frac{6}{p}+\frac{1}{p^2}\right)
\end{align*} and
\begin{align*}
V_{60}(B): =\int\limits_{\substack{(\boldsymbol{\eta}', \eta_8)\in\mathcal{R}(B)\\ \eta_6>|\eta_7|}}\frac{1}{\eta_5}\mathrm{d}\eta_1\cdots\mathrm{d}\eta_8.
\end{align*}
\end{lemma}
\begin{proof}
Firstly, we rewrite the result of Lemma \ref{L33} as follow.
\begin{align}\label{V32}
N_6(B) =\sum_{\substack{\boldsymbol{\eta}''' \in \mathbb{Z}_1\times\cdots\times\mathbb{Z}_5\times\mathbb{Z}_7}}
   \vartheta_6(\boldsymbol{\eta}''')V_6(\boldsymbol{\eta}''';B)+O\left(B(\log B)^4\right),
\end{align}
where\begin{align*}
\vartheta_6(\boldsymbol{\eta}'''):= \mathcal{A}(\vartheta'_1(\boldsymbol{\eta}'),\eta_6)\qquad\text{and}\qquad
V_6(\boldsymbol{\eta}''';B):= V_{16}(\boldsymbol{\eta}''';B).
\end{align*}

Similar to (\ref{E2}), by (\ref{Special2}) and (\ref{Quadratic1}), we imply easily
\begin{align}\label{E38}
\vartheta'_1(\boldsymbol{\eta}')=\prod\limits_{p}\vartheta'_{1,p}\left(I_p(\boldsymbol{\eta}')\right)
\end{align}
with $I_p(\boldsymbol{\eta}')=\left\{i\in\{1,\ldots,7\}:p\mid \eta_i\right\}$ and
\begin{align*}
\vartheta'_{1,p}(I)=
\begin{cases}
     1,  &I=\emptyset, \{1\}, \{5\}, \{7\},\\
    1-\frac{1}{p},&I=\{3\}, \{4\}, \{6\}, \{2,5\}, \{1,2\}, \{2,3\}, \{3,4\}, \{4,6\}, \{6,7\},\\
    1-\frac{2}{p}, &I =\{2\},\\
    0, & \text{all other} \;I \subset \{1, \dots, 7\}.
\end{cases}
\end{align*}
Thus, we have $\vartheta'_1(\boldsymbol{\eta}')\in\Theta'_{4,7}(2)$ \cite[Definition 7.8]{DE2}. By \cite[Fig.2]{DE2}, we get $\vartheta'_1(\boldsymbol{\eta}')\in\Theta_{2,7}(C)$ \cite[Definition 4.2]{DE2} for some $C\in\mathbb{Z}_{\geq0}$. Furthermore, we have
$\vartheta_6(\boldsymbol{\eta}''')= \mathcal{A}(\vartheta'_1(\boldsymbol{\eta}'),\eta_6)\in\Theta_{2,6}(C)$. Applying \cite[Lemma 5.1(3)]{DE2} to (\ref{H315}), we have
\begin{align*}
V_6(\boldsymbol{\eta}''';B)\ll\frac{B^{5/6}}{\eta_1^{1/2}\eta_3^{1/3}\eta_4^{2/3}\eta_5^{1/6}|\eta_7|^{7/6}}
=\frac{B}{\eta_1\cdots\eta_5|\eta_7|}\left(\frac{B}{\eta_1^3\eta_2^6\eta_3^4\eta_4^2\eta_5^5|\eta_7|^{-1}}\right)^{-1/6}.
\end{align*}
Futhermore, using (\ref{H311}) to bound $\eta_6$ and (\ref{H313}) to bound $\eta_8$, we have
\begin{align*}
V_6(\boldsymbol{\eta}''';B)\ll\frac{1}{\eta_5}\cdot\frac{B}{\eta_1^2\eta_2^4\eta_3^3\eta_4^2\eta_5^3}\cdot\frac{B}{\eta_1^2\eta_2^3\eta_3^2\eta_4\eta_5^2}
=\frac{B}{\eta_1\cdots\eta_5|\eta_7|}\left(\frac{B}{\eta_1^3\eta_2^6\eta_3^4\eta_4^2\eta_5^5|\eta_7|^{-1}}\right).
\end{align*}
Thus, we have
\begin{align*}
V_6(\boldsymbol{\eta}''';B)\ll\frac{B}{\eta_1\cdots\eta_5|\eta_7|}\min\left\{\left(\frac{B}{\eta_1^3\eta_2^6\eta_3^4\eta_4^2\eta_5^5|\eta_7|^{-1}}\right)^{-1/6},
\frac{B}{\eta_1^3\eta_2^6\eta_3^4\eta_4^2\eta_5^5|\eta_7|^{-1}}\right\}.
\end{align*}

Now, the conditions of \cite[Proposition 4.3]{DE2} are satisfied, and we apply it (see (\ref{A})) to the main term of (\ref{V32}) with $(r,s)=(5,1)$, we imply that
\begin{align*}
\sum_{\substack{\boldsymbol{\eta}''' \in \mathbb{Z}_1\times\cdots\times\mathbb{Z}_5\times\mathbb{Z}_7}}
  \vartheta_6(\boldsymbol{\eta}''')V_6(\boldsymbol{\eta}''';B)=\mathcal{A}( \vartheta_6(\boldsymbol{\eta}'''),\boldsymbol{\eta}''') &\int\limits_{\boldsymbol{\eta}'''}V_6(\boldsymbol{\eta}''';B)\mathrm{d}\boldsymbol{\eta}'''\\
 & +O\left(B(\log B)^4\log\log B\right).
\end{align*}
By (\ref{V31}), we have
\begin{align*}
\int\limits_{\boldsymbol{\eta}'''}V_6(\boldsymbol{\eta}''';B)\mathrm{d}\boldsymbol{\eta}'''=V_{60}(B).
\end{align*}

Since $\mathcal{A}( \vartheta_6(\boldsymbol{\eta}'''),\boldsymbol{\eta}''')=\mathcal{A}(\vartheta'_1(\boldsymbol{\eta}'),\boldsymbol{\eta}')$ and
$\vartheta'_1(\boldsymbol{\eta}')\in\Theta'_{4,7}(2)$, we get that the conditions of \cite[Corollary 7.10]{DE2} are satisfied.
Similar to (\ref{W}), by \cite[Corollary 7.10]{DE2} and (\ref{E38}), we have
\begin{align}\label{W0}
\mathcal{A}( \vartheta_6(\boldsymbol{\eta}'''),\boldsymbol{\eta}''')=\prod\limits_{p}\omega_p,
\end{align}
where
\begin{align*}
\omega_p=\left(1-\frac{1}{p}\right)^6\left(1+\frac{6}{p}+\frac{1}{p^2}\right).
\end{align*}
\end{proof}

\subsection{The asymptotic formula on $N_7(B)$}
\subsubsection{The second summation over $\eta_7$}
\begin{lemma}\label{L35}
We have
\begin{align*}
N_7(B)= \sum_{\substack{\boldsymbol{\eta} \in \mathbb{Z}_1\times\cdots\times\times\mathbb{Z}_6}}
   \mathcal{A}(\vartheta'_1(\boldsymbol{\eta}'),\eta_7)V_{17}(\boldsymbol{\eta};B)+O\left(B(\log B)^4\right),
\end{align*}
where $\mathcal{A}(\vartheta'_1(\boldsymbol{\eta}'),\eta_7)$ is defined in Section 2 and
\begin{align}\label{V33}
V_{17}(\boldsymbol{\eta};B) =\frac{1}{\eta_5}\int\limits_{\substack{(\boldsymbol{\eta}', \eta_8)\in\mathcal{R}(B)\\ \eta_6\leq|\eta_7|}}
\mathrm{d}\eta_7\mathrm{d}\eta_8.
\end{align}
\end{lemma}
\begin{proof}
This is similar to Lemma \ref{L33}. Firstly, by Lemma \ref{L32}, we have
\begin{align}\label{E39}
N_7(B)=\sum_{\substack{\boldsymbol{\eta} \in \mathbb{Z}_1\times\cdots\times\mathbb{Z}_6}}\sum_{\substack{k|\eta_2\\ \gcd(k,\eta_1\eta_3)}}
\frac{\mu(k)}{k}\Sigma,
\end{align}
where
\begin{align*}
\Sigma:=\sum_{\substack{\eta_7 \in \mathbb{Z}_7\\\eta_6\leq|\eta_7|}}\vartheta(\eta_7)\sum_{\substack{\varrho^2\equiv-\eta_3\eta_6\eta_7/\eta_1\pmod{k\eta_5}\\ 1\leq\varrho\leq k\eta_5\\ \gcd(\varrho,k\eta_5)}}g(\eta_7)
\end{align*}
with
\begin{align*}
\vartheta(\eta_7):= \widetilde{\vartheta}_1(\boldsymbol{\eta}',k) \qquad\text{and}\qquad g(\eta_7):= V_1(\boldsymbol{\eta},\eta_7;B).
\end{align*}

Similar to Lemma \ref{L33}, by (\ref{E32!}), (\ref{Special2}), we yield immediately that
\begin{align*}
\vartheta(\eta_7)\in \Theta_2(b,1,1,1)\qquad\text{and}\qquad b=\prod\limits_{p\mid \eta_1\eta_2\eta_3\eta_4\eta_5}p.
\end{align*}
However, by Lemma \ref{L11}, we can not get the desired error terms. The key observation is that, as in \cite[Section 5.4]{BR1}, we can replaced
$\vartheta_0(\boldsymbol{\eta}')$ (\ref{Special2}) by $\vartheta'_0(\boldsymbol{\eta}')$, where $\vartheta'_0(\boldsymbol{\eta}')=\prod\limits_{p}\vartheta'_{0,p}\left(I_p(\boldsymbol{\eta}')\right)$
with
\begin{align*}
\vartheta'_{0,p}(I)=
\begin{cases}
     1,  &\vartheta_{0,p}(I)=1\;\text{or}\;I=\{5,7\},\\
    0, & \text{all other} \;I \subset \{1, \dots, 7\}.
\end{cases}
\end{align*}
Note that we have been able to add the constraint $I=\{5,7\}$, since the sum over $\varrho$ is zero whenever $\gcd(\eta_5,\eta_7)>1$.
Now, by (\ref{E32!}) and $\vartheta'_0(\boldsymbol{\eta}')$, we yield immediately that
\begin{align*}
\vartheta(\eta_7)\in \Theta_2(b,1,1,1)\qquad\text{and}\qquad b=\prod\limits_{p\mid \eta_1\eta_2\eta_3\eta_4}p.
\end{align*}
Similar to (\ref{E34}),
\begin{align*}
g(\eta_7)\ll\left(\frac{B}{\eta_1\eta_5|\eta_7|}\right)^{1/2}.
\end{align*}
By Lemma \ref{L11}, we get that
\begin{align*}
\Sigma=\phi^\ast(k\eta_5)\mathcal{A}(\vartheta(\eta_7),\eta_7,k\eta_5)&\int_{\substack{\eta_7\in\mathbb{Z}_7\\ \eta_6\leq|\eta_7|}}g(\eta_7)\mathrm{d}\eta_7\\
&+O\left(kB^{1/2}\log B\frac{\eta_5^{1/2}2^{\omega(\eta_1\eta_2\eta_3\eta_4)}}{\eta_1^{1/2}\eta_6^{1/2}}\right).
\end{align*}
When summing the error term over $k$ and the remaining variables, we have
\begin{align*}
&\ll B^{1/2}\log B\sum_{\substack{\boldsymbol{\eta}\\(\ref{H311})}}
\frac{\eta_5^{1/2}2^{\omega(\eta_2)+\omega(\eta_1\eta_2\eta_3\eta_4)}}{\eta_1^{1/2}\eta_6^{1/2}}\\
&\ll B\log B\sum_{\substack{\eta_1,\cdots,\eta_4,\eta_6}}\frac{2^{\omega(\eta_2)+\omega(\eta_1\eta_2\eta_3\eta_4)}}{\eta_1^{3/2}\eta_2^2\eta_3^{3/2}\eta_4\eta_6}\\
&\ll B(\log B)^4.
\end{align*}

Now, we have
\begin{align*}
N_7(B) = \sum_{\substack{\boldsymbol{\eta} \in \mathbb{Z}_1\times\cdots\times\mathbb{Z}_6}}\sum_{\substack{k|\eta_2\\ \gcd(k,\eta_1\eta_3)}}
\frac{\mu(k)}{k}\phi^\ast(k\eta_5)&\mathcal{A}(\vartheta(\eta_7),\eta_7,k\eta_5)V_{17}(\boldsymbol{\eta};B)\\
&+O\left(B(\log B)^4\right).
\end{align*}
In order to complete this lemma, it suffices to prove
\begin{align*}
\sum_{\substack{k|\eta_2\\ \gcd(k,\eta_1\eta_3)}}\frac{\mu(k)}{k}\phi^\ast(k\eta_5)&\mathcal{A}(\vartheta(\eta_7),\eta_7,k\eta_5)
= \mathcal{A}(\vartheta'_1(\boldsymbol{\eta}'),\eta_7).
\end{align*}
Similar to (\ref{E36}), we can prove it easily.
\end{proof}

\subsubsection{The remaining summations}
\begin{lemma}\label{L36}
We have
\begin{align*}
N_7(B)=\left(\prod\limits_{p}\omega_p\right)V_{70}(B)+O\left(B(\log B)^4\log\log B\right),
\end{align*}
where
\begin{align*}
\omega_p:=\left(1-\frac{1}{p}\right)^6\left(1+\frac{6}{p}+\frac{1}{p^2}\right)
\end{align*} and
\begin{align*}
V_{70}(B): =\int\limits_{\substack{(\boldsymbol{\eta}', \eta_8)\in\mathcal{R}(B)\\ \eta_6\leq|\eta_7|}}\frac{1}{\eta_5}\mathrm{d}\eta_1\cdots\mathrm{d}\eta_8.
\end{align*}
\end{lemma}
\begin{proof}
This is similar to Lemma \ref{L34}. Firstly, we rewrite the result of Lemma \ref{L35} as follow.
\begin{align}\label{V34}
N_7(B)=\sum_{\substack{\boldsymbol{\eta}\in \mathbb{Z}_1\times\cdots\times\mathbb{Z}_6}}
   \vartheta_7(\boldsymbol{\eta})V_7(\boldsymbol{\eta};B)+O\left(B(\log B)^4\right),
\end{align}
where\begin{align*}
\vartheta_7(\boldsymbol{\eta}):= \mathcal{A}(\vartheta'_1(\boldsymbol{\eta}'),\eta_7)\qquad\text{and}\qquad
V_7(\boldsymbol{\eta};B):= V_{17}(\boldsymbol{\eta};B).
\end{align*}

From (\ref{E38}), we have $\vartheta'_1(\boldsymbol{\eta}')\in\Theta'_{4,7}(2)$ \cite[Definition 7.8]{DE2}. By \cite[Fig.2]{DE2}, we get $\vartheta'_1(\boldsymbol{\eta}')\in\Theta_{2,7}(C)$ \cite[Definition 4.2]{DE2} for some $C\in\mathbb{Z}_{\geq0}$. Furthermore, we have
$\vartheta_7(\boldsymbol{\eta})= \mathcal{A}(\vartheta'_1(\boldsymbol{\eta}'),\eta_7)\in\Theta_{2,6}(C)$. Applying \cite[Lemma 5.1(6)]{DE2} to (\ref{H315}), we have
\begin{align*}
V_7(\boldsymbol{\eta};B)\ll\frac{B^{3/4}}{\eta_1^{1/2}\eta_3^{1/4}\eta_4^{1/2}\eta_5^{1/4}\eta_6^{3/4}}
=\frac{B}{\eta_1\cdots\eta_6}\left(\frac{B}{\boldsymbol{\eta}^{(2,4,3,2,3,1)}}\right)^{-1/4}.
\end{align*}

Now, the conditions of \cite[Proposition 4.3]{DE2} are satisfied, and we apply it (see (\ref{A})) to the main term of (\ref{V34}) with $(r,s)=(5,1)$, we imply that
\begin{align*}
\sum_{\substack{\boldsymbol{\eta} \in \mathbb{Z}_1\times\cdots\times\mathbb{Z}_6}}
  \vartheta_7(\boldsymbol{\eta})V_7(\boldsymbol{\eta};B)=\mathcal{A}( \vartheta_7(\boldsymbol{\eta}),\boldsymbol{\eta}) &\int\limits_{\boldsymbol{\eta}}V_7(\boldsymbol{\eta};B)\mathrm{d}\boldsymbol{\eta}
 +O\left(B(\log B)^4\log\log B\right).
\end{align*}
By (\ref{V33}), we have
\begin{align*}
\int\limits_{\boldsymbol{\eta}}V_7(\boldsymbol{\eta};B)\mathrm{d}\boldsymbol{\eta}=V_{70}(B).
\end{align*}
Similar to (\ref{W0}), we can check that the conditions of \cite[Corollary 7.10]{DE2} are satisfied.  By (\ref{E38}), we have
\begin{align*}
\mathcal{A}( \vartheta_7(\boldsymbol{\eta}),\boldsymbol{\eta})=\prod\limits_{p}\omega_p,
\end{align*}
where
\begin{align*}
\omega_p=\left(1-\frac{1}{p}\right)^6\left(1+\frac{6}{p}+\frac{1}{p^2}\right).
\end{align*}
\end{proof}

\subsection{The final proof}
\begin{lemma}\label{L37}
We have
\begin{align*}
N_{U_3,H}(B)=\left(\prod\limits_{p}\omega_p\right)V_0(B)+O\left(B(\log B)^4\log\log B\right),
\end{align*}
where
\begin{align*}
\omega_p:=\left(1-\frac{1}{p}\right)^6\left(1+\frac{6}{p}+\frac{1}{p^2}\right)
\end{align*} and
\begin{align*}
V_0(B): =\int\limits_{\substack{(\boldsymbol{\eta}', \eta_8)\in\mathcal{R}(B)}}\frac{1}{\eta_5}\mathrm{d}\eta_1\cdots\mathrm{d}\eta_8.
\end{align*}
\end{lemma}
\begin{proof}
This follows from (\ref{N}), Lemma \ref{L34} and Lemma \ref{L36}.
\end{proof}

By \cite[Theorem 4]{DE4} and \cite[Theorem 1.3]{DEJOTE1}, it is easy to compute the constant $\alpha(\widetilde{S}_3)$. We have
\begin{align}\label{W3}
\alpha(\widetilde{S}_3)=\frac{1}{180}\cdot\frac{1}{\#W(\mathbf{A}_4)}=\frac{1}{21600},
\end{align}
where we have used $\#W(\mathbf{A}_n)=(n + 1)!$. Let
\begin{equation*}
\omega_\infty(\widetilde{S}_3):=3\int_{\max\{|z_0^3|,
|z_0z_2z_3|, |z_0^2z_2|, |z_0^2z_3|, |z_3(z_2^2+z_0z_3)|\}\leq 1}\mathrm{d}z_0 \mathrm{d}z_1 \mathrm{d}z_2.
\end{equation*}
For $(\eta_1,\eta_2,\eta_4,\eta_5,\eta_6)\in\mathbb{R}^5_{\geq 1}$, we introduce the conditions
\begin{align}
&\eta_1^2\eta_2^4\eta_4^2\eta_5^3\eta_6\leq B,\label{H316}\\
&\eta_1^2\eta_2^4\eta_4^2\eta_5^3\eta_6\leq B,\qquad \eta_1^{-1}\eta_2^{-2}\eta_4^2\eta_5^{-3}\eta_6^4\leq B.\label{H317}
\end{align}

\begin{lemma}\label{L38}
 Let $\alpha(\widetilde{S}_3)$ be as in (\ref{W3}), we have
  \begin{equation*}
    V_0'(B) := \int_{\substack{(\boldsymbol{\eta}', \eta_8)\in \mathcal{R}'(B)\\
    (\ref{H317})}} \frac{1}{\eta_5} \mathrm{d}\eta_1\cdots\mathrm{d}\eta_8,
  \end{equation*}
where\begin{align*}
\mathcal{R}'(B)=\{(\boldsymbol{\eta}', \eta_8)\in J_1'\times\cdots\times J_8'\mid(\ref{H311})-(\ref{H315})\}.
\end{align*}
Then $\frac{1}{21600}\omega_\infty(\widetilde{S}_3) B(\log B)^5 = V_0'(B)$.
\end{lemma}
\begin{proof}
From (\ref{P31}), we have
\begin{equation}\label{Y2}
\begin{aligned}
(y_0: y_1: y_2)&\mapsto\\
&\left(y_0^3: y_0y_1y_2: y_0^2y_1: y_0^2y_2: -y_2(y_1^2+y_0y_2)\right).
\end{aligned}
\end{equation}
Inserting $y_0 = \eta_1\eta_2^2\eta_3\eta_4\eta_5^2\eta_6$, $y_1 =\eta_1\eta_2\eta_5\eta_8$, $y_2= \eta_4\eta_6^2\eta_7$, $-(y_1^2+y_0y_2) =
\eta_1\eta_2^2\eta_5^3\eta_9$ into (\ref{Y2}) and cancelling out $\eta_1\eta_2^2\eta_4\eta_5^3\eta_6^2$ gives $\phi''$ as
in the proof of Lemma \ref{L31}. Let $\eta_1, \eta_2,  \eta_4, \eta_5, \eta_6\in\mathbb{R}_{\neq 0}$ and $\eta_3, \eta_7, \eta_8\in\mathbb{R}$. With $l=B\eta_1\eta_2^2\eta_4\eta_5^3\eta_6^2$, we utilize the coordinate transformation
\begin{align*}
z_0=l^{-1/3}\eta_1\eta_2^2\eta_4\eta_5^2\eta_6\cdot\eta_3,\qquad z_1=l^{-1/3}\eta_1\eta_2\eta_5\cdot\eta_8,\qquad z_2=l^{-1/3}\eta_4\eta_6^2\cdot\eta_7
\end{align*}
of Jacobi determinant
\begin{align*}
\frac{\eta_1\eta_2\eta_4\eta_5\eta_6}{B}\cdot\frac{1}{\eta_5}
\end{align*}
and get
\begin{align}\label{W4}
\omega_\infty(\widetilde{S}_3)=\frac{3\eta_1\eta_2\eta_4\eta_5\eta_6}{B}\int_{(\boldsymbol{\eta}', \eta_8)\in \mathcal{R}''(B)}\frac{1}{\eta_5}\mathrm{d}\eta_3 \mathrm{d}\eta_7\mathrm{d}\eta_8,
\end{align}
where
\begin{align*}
\mathcal{R}''(B)=\{(\boldsymbol{\eta}', \eta_8)\in\mathbb{R}^2_{\neq 0}\times\mathbb{R}\times\mathbb{R}^3_{\neq 0}\times\mathbb{R}^2\mid(\ref{H311})-(\ref{H315})\}.
\end{align*}

An application of Lemma \ref{L12} with exchanged roles of $\eta_3$ and $\eta_6$ gives
\begin{align}\label{a3}
\alpha(\widetilde{S}_3)(\log B)^5=\frac{1}{3}\int\limits_{\substack{(\eta_1, \eta_2,  \eta_4, \eta_5, \eta_6)\in\mathbb{R}^5_{\geq 1}\\ (\ref{H317})}}
\frac{\mathrm{d}\eta_1\mathrm{d}\eta_2\mathrm{d}\eta_4\mathrm{d}\eta_5\mathrm{d}\eta_6}{\eta_1\eta_2\eta_4\eta_5\eta_6},
\end{align}
since $[-K_{\widetilde{S}_3}] = [2E_1+4E_2+3E_3+2E_4+3E_5+E_6]$, $[E_7]=[E_1+2E_2+E_3+2E_5-E_6]$. From (\ref{W3}) and (\ref{W4})-(\ref{a3}), we have
\begin{align*}
\frac{1}{21600}\omega_\infty(\widetilde{S}_3) B(\log B)^5 = V_0'(B).
\end{align*}
\end{proof}

To finish our proof, we compare $V_0(B)$ from Lemma \ref{L37} with $V_0'(B)$ from Lemma \ref{L38}.
\begin{lemma}\label{L39}
 We have
\begin{equation*}
V_0(B) = V_0'(B)+O\left(B(\log B)^4\right).
\end{equation*}
\end{lemma}
\begin{proof}
 Let
\begin{align*}
  \mathcal{D}_0(B) &:= \{(\boldsymbol{\eta}', \eta_8) \in \mathbb{R}^8 \, \text{and}\,(\ref{H311})-(\ref{H315})
  \mid \eta_1,\ldots,\eta_6,|\eta_7| \geq 1\},\\
  \mathcal{D}_1(B) &:= \{(\boldsymbol{\eta}', \eta_8) \in \mathbb{R}^8 \, \text{and}\,(\ref{H311})-(\ref{H315})
  \mid \eta_1,\ldots,\eta_6,|\eta_7| \geq 1, (\ref{H316})\},\\
   \mathcal{D}_2(B) &:= \{(\boldsymbol{\eta}', \eta_8) \in \mathbb{R}^8 \, \text{and}\,(\ref{H311})-(\ref{H315})
  \mid \eta_1,\ldots,\eta_6,|\eta_7| \geq 1, (\ref{H317})\},\\
   \mathcal{D}_3(B) &:= \{(\boldsymbol{\eta}', \eta_8) \in \mathbb{R}^8 \, \text{and}\,(\ref{H311})-(\ref{H315})
  \mid \eta_1,\ldots,\eta_6\geq 1, (\ref{H317})\},\\
 \mathcal{D}_3(B) &:= \{(\boldsymbol{\eta}', \eta_8) \in \mathbb{R}^8 \, \text{and}\,(\ref{H311})-(\ref{H315})
  \mid \eta_1,\eta_2,\eta_4,\eta_5,\eta_6\geq 1, (\ref{H317})\}.
\end{align*}
For $i\in\{0,\ldots,4\}$, let
\begin{align*}
V_i(B):=\int_{(\boldsymbol{\eta}', \eta_8)\in\mathcal{D}_i(B)}\frac{\mathrm{d}\eta_1\cdots\mathrm{d}\eta_8}{\eta_5}.
\end{align*}
Then $V_0(B)$ is as in Lemma \ref{L37} and $V_4(B) = V_0'(B)$. It suffices to show that $V_i(B)-V_{i-1}(B) = O(B(\log B)^4)$ for $1 \leq i \leq 4$. This holds for $i = 1$, since, by (\ref{H311}) and $\eta_3 \geq 1$, we have $\mathcal{D}_1(B) =\mathcal{D}_0(B)$.

Moreover, using  \cite[Lemma 5.1(6)]{DE2} and (\ref{H315}) to bound the integral over $\eta_7$ and $\eta_8$, we have
\begin{align*}
V_2(B) - V_1(B)& \ll \int\limits_{\substack{\eta_1,\ldots,\eta_6\geq 1 \\ \eta_1\eta_2^2\eta_3^2\eta_4^2\eta_5\eta_6^2\leq B \\ \eta_1^{-1}\eta_2^{-2}\eta_4^2\eta_5^{-3}\eta_6^4> B }}\frac{B^{3/4}}{\eta_1^{1/2}\eta_3^{1/4}\eta_4^{1/2}\eta_5^{1/4}\eta_6^{3/4}}\mathrm{d} \eta_1 \cdots \mathrm{d} \eta_6 \\
&\ll \int\limits_{\substack{\eta_1,\eta_2,\eta_4,\eta_5,\eta_6\geq 1\\ \eta_1^{-1}\eta_2^{-2}\eta_4^2\eta_5^{-3}\eta_6^4> B}}
\frac{B^{9/8}}{\eta_1^{7/8}\eta_2^{3/4}\eta_4^{5/4}\eta_5^{5/8}\eta_6^{3/2}}
\mathrm{d} \eta_1 \mathrm{d} \eta_2\mathrm{d} \eta_4\mathrm{d} \eta_5 \mathrm{d} \eta_6 \\
&\ll \int\limits_{\substack{\eta_2,\eta_4,\eta_5,\eta_6\geq 1}}
\frac{B}{\eta_2\eta_4\eta_5\eta_6}\mathrm{d} \eta_2\mathrm{d} \eta_4\mathrm{d} \eta_5 \mathrm{d} \eta_6 \\
&\ll B(\log B)^4,
\end{align*}
where $ \eta_1\eta_2^2\eta_3^2\eta_4^2\eta_5\eta_6^2\leq B$ comes from (\ref{H314}).

Using  \cite[Lemma 5.1(1)]{DE2} and (\ref{H315}) to bound the integral over $\eta_8$, we have
\begin{align*}
V_3(B) - V_2(B)& \ll \int\limits_{\substack{\eta_1,\ldots,\eta_6\geq 1 \\ |\eta_7|<1,(\ref{H311}),(\ref{H317})}}\frac{B^{1/2}}{(\eta_1\eta_5|\eta_7|)^{1/2}}\mathrm{d} \eta_1 \cdots \mathrm{d} \eta_7 \\
&\ll \int\limits_{\substack{\eta_1,\eta_2,\eta_4,\eta_5,\eta_6\geq 1 \\ (\ref{H317})}}\frac{B^{5/6}}{\eta_1^{7/6}\eta_2^{4/3}\eta_4^{2/3}\eta_5^{3/2}\eta_6^{1/3}}
\mathrm{d}\eta_1\mathrm{d}\eta_2\mathrm{d}\eta_4\mathrm{d}\eta_5\mathrm{d}\eta_6\\
&\ll \int\limits_{\substack{\eta_1,\eta_2,\eta_4,\eta_5\geq 1}}
\frac{B}{\eta_1\eta_2\eta_4\eta_5}\mathrm{d}\eta_1\mathrm{d}\eta_2\mathrm{d}\eta_4\mathrm{d}\eta_5\\
&\ll B(\log B)^4.
\end{align*}

Finally, using  \cite[Lemma 5.1(6)]{DE2} and (\ref{H315}) to bound the integral over $\eta_7$ and $\eta_8$, we have
\begin{align*}
V_4(B)-V_3(B)& \ll \int\limits_{\substack{\eta_1,\eta_2,\eta_4,\eta_5,\eta_6\geq 1\\ \eta_3<1,(\ref{H316})}}\frac{B^{3/4}}{\eta_1^{1/2}\eta_3^{1/4}\eta_4^{1/2}\eta_5^{1/4}\eta_6^{3/4}}\mathrm{d} \eta_1 \cdots \mathrm{d} \eta_6 \\
&\ll \int\limits_{\substack{\eta_1,\eta_2,\eta_4,\eta_5\geq 1}}\frac{B}{\eta_1\eta_2\eta_4\eta_5}
\mathrm{d}\eta_1\mathrm{d}\eta_2\mathrm{d}\eta_4\mathrm{d}\eta_5\\
&\ll B(\log B)^4.
\end{align*}
\end{proof}

Theorem 1.3 follows from Lemma \ref{L37}, Lemma \ref{L38} and Lemma \ref{L39}.

\subsection*{Acknowledgements}
The author would like to thank Professor Valentin Blomer and Guangshi L\"{u} for their constant encouragement and helpful comments. This work was completed when the author stayed in University of Bonn. The author sincerely expresses his deep gratitude to Professor Valentin Blomer for the support and generous hospitality. Moreover, the author is grateful to Professor Ulrich Derenthal for the comments and suggestons on an earlier draft of this paper. It is also a pleasure to thank Professor R\'{e}gis de la Bret\`{e}che for suggesting the  modifications of Table \ref{table0}.
This work is supported by the China Scholarship Council.

\bibliographystyle{amsplain}

\end{document}